\documentclass[a4paper,12pt]{article}
\usepackage{amsmath}
\usepackage{amsfonts}
\usepackage{amssymb}
\usepackage{latexsym}
\usepackage{epsfig}
\usepackage{graphicx}
\usepackage{oldgerm}
\usepackage{theorem}
\def\C{\mathbb{C}}
\def\E{{\mathbb E}}
\def\mbK{\mathbb{K}}
\def\N{\mathbb{N}}
\def\P{{\mathbb P}}

\def\R{\mathbb{R}}

\def\X{\mib{X}}
\def\Y{\mib{Y}}
\def\Z{\mathbb{Z}}

\def\f{\mib{f}}

\def\t{\mib{t}}

\def\x{\mib{x}}
          \def\X{\mib{X}}

\def\0{\mib{0}}
\def\rC{{\rm C}}

\def\mM{\mathfrak{M}}

\def\mY{\mathfrak{Y}}

\def\bK{{\bf K}}
\def\Ai{{\rm Ai}}

\def\cA{{\cal A}}

\def\1{{\bf 1}}
\def\fd{\stackrel{\rm f.d.}{\longrightarrow}}
\def\Det{{\rm Det}}
\def\supp{{\rm supp}\ }
\theorembodyfont{\itshape}
\newtheorem{thm}{Theorem}[section]
\newtheorem{lem}[thm]{Lemma}

\newtheorem{prop}[thm]{Proposition}

\newtheorem{df}[thm]{Definition}
\newcommand{\mib}[1]{\mbox{\boldmath $#1$}}
\newcommand{\SSC}[1]{\section{#1}\setcounter{equation}{0}}
\newcommand{\qed}{\hbox{\rule[-2pt]{3pt}{6pt}}}

\begin{document}

\title{Markov property of determinantal processes \\
with extended sine, Airy, and Bessel kernels}
\author{
Makoto Katori
\footnote{
Department of Physics,
Faculty of Science and Engineering,
Chuo University, 
Kasuga, Bunkyo-ku, Tokyo 112-8551, Japan;
e-mail: katori@phys.chuo-u.ac.jp
}
and 
Hideki Tanemura
\footnote{
Department of Mathematics and Informatics,
Faculty of Science, Chiba University, 
1-33 Yayoi-cho, Inage-ku, Chiba 263-8522, Japan;
e-mail: tanemura@math.s.chiba-u.ac.jp
}}
\date{22 June 2011}
\pagestyle{plain}
\maketitle
\begin{abstract}
When the number of particles is finite, 
the noncolliding Brownian motion (the Dyson model) and
the noncolliding squared Bessel process are determinantal
diffusion processes for any deterministic initial
configuration $\xi=\sum_{j \in \Lambda} \delta_{x_j}$, 
in the sense that any multitime correlation function
is given by a determinant associated with the
correlation kernel, which is specified by an
entire function $\Phi$ having zeros in $\supp \xi$.
Using such entire functions $\Phi$,
we define new topologies called the 
$\Phi$-moderate topologies.
Then we construct three infinite-dimensional determinantal
processes, as the limits of sequences of determinantal 
diffusion processes with finite numbers of particles
in the sense of finite dimensional distributions
in the $\Phi$-moderate topologies, so that
the probability distributions are continuous
with respect to initial configurations $\xi$
with $\xi(\R)=\infty$.
We show that our three infinite particle systems are versions of
the determinantal processes with the extended sine,
Bessel, and Airy kernels, respectively,
which are reversible with respect to the
determinantal point processes obtained in
the bulk scaling limit and the soft-edge scaling limit of
the eigenvalue distributions of the Gaussian unitary
ensemble, and the hard-edge scaling limit of that
of the chiral Gaussian unitary ensemble studied in
the random matrix theory.
Then Markovianity is proved for the three infinite-dimensional
determinantal processes. 

\vskip 0.3cm
\noindent{\bf Keywords} \,
Determinantal processes $\cdot$ Correlation kernels $\cdot$
Random matrix theory $\cdot$ Infinite particle systems $\cdot$
Markov property $\cdot$ Entire function and topology

\vskip 0.5cm
\noindent{\bf Mathematics Subject Classification (2010)} \,
15B52, 30C15, 47D07, 60G55, 82C22

\end{abstract}

\clearpage

\SSC{Introduction}\label{chap: Introduction}

Let $\mM$ be the space of nonnegative 
integer-valued Radon measures on $\R$,
which is a Polish space with the {\it vague topology}:
let ${\rm C}_0(\R)$ be the set of all 
continuous real-valued functions with compact supports, and 
we say $\xi_n, n\in\N$ converges to $\xi$ vaguely, if 
$\lim_{n \to \infty} \int_{\R} \varphi(x) \xi_n(dx)
=\int_{\R} \varphi(x) \xi(dx)$ 
for any $\varphi \in {\rm C}_0(\R)$.
Each element $\xi$ of $\mM$ can be represented as
$\xi(\cdot) = \sum_{j\in \Lambda}\delta_{x_j}(\cdot)$
with an index set $\Lambda$ and a sequence of points 
$\x =(x_j)_{j \in \Lambda}$ in $\R$ 
satisfying $\xi(K)=\sharp\{x_j: x_j \in K\} < \infty$ 
for any compact subset $K \subset \R$.
We call an element $\xi$ of $\mM$ an unlabeled configuration,
and a sequence of points $\x$ a labeled configuration.
A probability measure on a configuration space is called 
a determinantal point process or Fermion point process,
if its correlation functions are generally represented 
by determinants \cite{ST03, Sos00}.
In the present paper we say that an $\mM$-valued process $\Xi(t)$ 
is {\it determinantal}, 
if the multitime correlation functions 
for any chosen series of times are represented by determinants. 
In other words, a {\it determinantal process} is 
an $\mM$-valued process such that, 
for any integer $M\in \N=\{1,2,\dots\}$,
$\f=(f_1,f_2,\dots,f_M) \in \rC_0(\R)^M$, a sequence of times
$\t=(t_1,t_2,\dots,t_M)$ with $0 < t_1 < \cdots < t_M < \infty$,
if we set $\chi_{t_m}(x)=e^{f_{m}(x)}-1, 1 \leq m \leq M$,
the moment generating function of multitime distribution, 
$
{\Psi}^{\t}[\f]
\equiv \E \left[\exp \left\{ \sum_{m=1}^{M} 
\int_{\R} f_m(x) \Xi(t_m, dx) \right\} \right],
$
is given by a Fredholm determinant
\begin{equation}
{\Psi}^{\t}[\f]
= \mathop{\Det}_
{\substack
{(s,t)\in \{t_1, t_2, \dots, t_M \}^2, \\
(x,y)\in \R^2}
}
=\Det \Big[\delta_{s t} \delta(x-y)
+ \mbK(s, x; t, y) \chi_{t}(y) \Big],
\label{eqn:Fred}
\end{equation}
with a locally integrable function $\mbK$ called a 
{\it correlation kernel} \cite{KT07b, KT10a}.
Finite- and infinite-dimensional determinantal processes
were introduced as multi-matrix models \cite{EM98,NF98}, 
tiling models \cite{Joh02}, 
and surface growth models \cite{PS02},
and they have been extensively studied.
In the present paper we study the infinite-dimensional determinantal
processes describing one-dimensional infinite particle systems
with long-range repulsive interactions.
They are obtained by taking appropriate $N \to \infty$
limits of the $N$-particle systems of
{\it noncolliding diffusion processes},
which dynamically simulate the eigenvalue statistics
of the Gaussian random matrix ensembles studied in
the random matrix theory \cite{Meh04,For10}.
The purpose of the present paper is to prove
{\it Markov property} of the three infinite-dimensional
determinantal processes with the correlation kernels
called the {\it extended sine, Airy, and Bessel kernels},
which are reversible with respect to the probability measures
obtained in the bulk scaling limit and the soft-edge scaling limit
of the eigenvalue distribution in the 
{\it Gaussian unitary ensemble} (GUE),
and in the hard-edge scaling limit of that in the
{\it chiral Gaussian unitary ensemble} (chGUE),
respectively. See, for instance, \cite{KT11b}
and the references therein.

Dyson \cite{Dys62} introduced a stochastic model of particles in $\R$
with a log-potential,
which obeys the stochastic differential equations (SDEs):
\begin{equation}
dX_j(t)=dB_j(t) + \sum_{1 \leq k \leq N, k \not= j}\frac{dt}{X_j(t)-X_k(t)}, 
\quad 1 \leq j \leq N, \quad t \in [0, \infty),
\label{eqn:Dyson}
\end{equation}
where $B_j(t)$'s are independent one-dimensional
standard Brownian motions.
This is a special case with the parameter $\beta=2$ of
Dyson's Brownian motion models \cite{Dys62,Meh04} but we will call it 
simply the Dyson model in this paper.
It is equivalent to a system of one-dimensional Brownian motions 
conditioned never to collide with each other \cite{KT11b}. 
For the solution $X_j(t)$, $j=1,2,\dots, N$ of (\ref{eqn:Dyson})
with initial values $X_j(0)=x_j$, $j=1,2,\dots, N$, 
we put $\Xi(t, \cdot) = \sum_{j=1}^N \delta_{X_j(t)}(\cdot)$ 
and $\xi^N (\cdot)= \sum_{j=1}^N \delta_{x_j}(\cdot)$. 
The (unlabeled) Dyson model starting from $\xi^{N}$
is denoted by $(\{\Xi(t)\}_{ t\in [0,\infty)},\P^{\xi^N})$.
In the previous paper \cite{KT10a} 
we showed that for any fixed unlabeled configuration $\xi^N$, 
the process $(\{\Xi(t)\}_{ t\in [0,\infty)}, \P^{\xi^N})$ 
is determinantal with the correlation kernel 
$\mbK^{\xi^{N}}$ specified by an entire function $\Pi_0$, which is
explicitly given by (\ref{def:K_xi}) below in the present paper.
Then for configurations $\xi \in \mM$ with $\xi(\R)=\infty$,
a sufficient condition was given so that 
the sequence of the processes 
$(\{\Xi(t)\}_{ t\in [0,\infty)},\P^{\xi\cap [-L,L]})$ 
converges to an $\mM$-valued determinantal process 
as $L\to\infty$ in the sense of finite dimensional distributions 
in the vague topology, where $\xi \cap [-L,L]$ denotes 
the restricted measure of $\xi$ on an interval $[-L,L]$.
Note that 
$\xi\in \mM$ implies $\xi\cap [-L, L](\R) <\infty$
for $0< L < \infty$.
The limit process is 
{\it the Dyson model with an infinite number of particles} 
and denoted by $(\{\Xi(t)\}_{ t\in [0,\infty)},\P^\xi)$. 
Recently, in \cite{KT10b} the tightness of the 
sequence of processes was proved. 
The class of configurations satisfying our condition, denoted by $\mY$, 
is large enough to carry the Poisson point processes, 
Gibbs states with regular conditions,
as well as the determinantal (Fermion) point process $\mu_{\sin}$
with the sine kernel:
\begin{equation}
K_{\sin}(x,y)\equiv
\frac{1}{2 \pi} \int_{|k| \leq \pi} dk \,
e^{i k(y-x)}
= \frac{\sin \{ \pi(y-x) \} }{\pi (y-x)},
\quad x, y \in \R,
\label{def:sine-kernel}
\end{equation}
where $i=\sqrt{-1}$.

From the uniqueness of solutions of (\ref{eqn:Dyson}), 
the Dyson model with a finite number of particles 
is a diffusion process ({\it i.e.} it is a strong Markov process
having a continuous path almost surely). 
Although the process $(\{\Xi(t)\}_{ t\in [0,\infty)}, \P^{\xi})$ 
is given as the limit of a sequence of diffusion processes 
$(\{\Xi(t)\}_{ t\in [0,\infty)}, \P^{\xi\cap [-L,L]})$, $L\in\N$, 
the Markov property could be lost in it.
In the Dyson model, because of the long-range interaction, 
if the number of particles is infinite, 
the probability distribution $\P^{\xi}$ is
not vaguely continuous with respect to the initial configuration $\xi$. 
Therefore the Markovianity of the process is not readily 
concluded in the infinite-particle limit of 
a sequence of Markov processes in the vague topology.

For a probability measure $\mu$ on $\mM$ 
we put $\P_{\mu}(\cdot)= \int_{\mM}\mu(d\xi)\P^{\xi}(\cdot)$. 
Suppose that $\xi(\R)=\infty$, $\mu$-almost surely.
For proving Markov property of the process 
$(\{\Xi(t)\}_{ t\in [0,\infty)}, \P_{\mu})$, 
it is sufficient to find a subset $\widehat{\mY}$ of $\mM$ 
with a topology ${\cal T}$, 
which is stronger than the vague topology, 
and a sequence of probability measures $\{\mu_N \}_{N\in\N}$ such that
\begin{enumerate}
\item $\P^{\xi}(\cdot)$ is continuous 
with respect to $\xi$ in the topology ${\cal T}$,

\item $\P_{\mu}(\Xi(t) \in \widehat{\mY})=1$, for any $t\in [0,\infty)$, 

\item $\mu_N (\{\eta \in \widehat{\mY} : \eta(\R)=N \})=1$, $N\in\N$, 

\item $(\{\Xi(t)\}_{ t\in [0,\infty)}, \P_{\mu_N})$ 
converges $(\{\Xi(t)\}_{ t\in [0,\infty)}, \P_{\mu})$ 
as $N\to\infty$ in the sense of finite dimensional distributions 
in the topology ${\cal T}$.
\end{enumerate}
In this paper first we will show that these conditions are 
satisfied in the case that $\mu=\mu_{\sin}$, if we use 
$\mY$ as $\widehat{\mY}$ with the topology ${\cal T}$ 
called the {\it $\Phi_0$-moderate topology} 
defined by (\ref{eqn:THC2}) given below.
We also show that the process 
$(\{\Xi(t)\}_{ t\in [0,\infty)}, \P_{\mu_{\sin}})$ 
is a {\it version} of the determinantal process 
$(\{\Xi(t)\}_{t \in [0, \infty)},{\bf P}_{\sin})$ 
associated with the {\it extended sine kernel} 
${\bf K}_{\sin}$ with density 1:
\begin{eqnarray}
{\bf K}_{\sin}(s,x;t,y) &\equiv&
\frac{1}{2 \pi} \int_{|k| \leq \pi} dk \,
e^{k^2(t-s)/2 + i k (y-x)}- {\bf 1}(s>t) p(s-t, x|y) 
\nonumber\\
&=& 
\left\{ \begin{array}{ll} 
\displaystyle{
\int_{0}^{1} du \, e^{\pi^2 u^2 (t-s)/2} 
\cos \{ \pi u (y-x)\} }
& \mbox{if $t>s $} \cr
K_{\sin}(x,y)
& \mbox{if $t=s$} \cr
\displaystyle{
- \int_{1}^{\infty} du \, 
e^{\pi^2 u^2 (t-s)/2} \cos \{ \pi u (y-x) \} }
& \mbox{if $t<s$}.
\end{array} \right.
\label{def:ex_sine-kernel}
\end{eqnarray}
(Two processes having the same state space are said to be
{\it equivalent} if they have the same finite-dimensional distributions,
and we also say that each one is a {\it version} of other
or they are versions of the same process.
See Section 1.1 of \cite{RY98}.)
Hence, we conclude that the determinantal process 
with the extended sine kernel is a Markov process 
which is reversible with respect to
the determinantal point process $\mu_{\sin}$ 
(Theorem \ref{Theorem:Dyson}).

The above strategy will be also used to 
show Markovianity of the infinite-dimensional determinantal process 
$(\{\Xi(t)\}_{ t\in [0,\infty)},{\bf P}_{J_\nu})$ associated 
with the {\it extended Bessel kernel} $\bK_{J_{\nu}}$ :
\begin{equation}
\bK_{J_{\nu}}(s,x;t,y) = \left\{
   \begin{array}{ll}
\displaystyle{
\int_{0}^{1} d u \,
e^{-2u(s-t)} J_{\nu}(2 \sqrt{u x})
J_{\nu}(2 \sqrt{u y})
} 
& \mbox{if} \quad s < t  \\
& \\
\displaystyle{
\frac{J_{\nu}(2 \sqrt{x}) \sqrt{y} J_{\nu}'(2 \sqrt{y})
-\sqrt{x} J_{\nu}'(2\sqrt{x}) J_{\nu}(2\sqrt{y})}{x-y}
}
& \mbox{if} \quad t=s \\
& \\
\displaystyle{
- \int_{1}^{\infty} d u \,
e^{-2u(s-t)} J_{\nu}(2 \sqrt{ux})
J_{\nu}(2 \sqrt{uy})
}
& \mbox{if} \quad s > t,
   \end{array} \right. 
\label{eqn:KBessel1}
\end{equation}
$x, y \in \R_+ \equiv \{x \in \R : x \geq 0 \}$, 
where $J_{\nu}(\cdot)$ is the Bessel function with index $\nu>-1$ defined by
\begin{equation}
J_{\nu}(z) = \sum_{n=0}^{\infty} 
\frac{(-1)^n}{\Gamma(n+1) \Gamma(n+1+\nu)}
\left( \frac{z}{2} \right)^{2n+\nu},
\quad z\in \C
\nonumber
\end{equation}
with the gamma function 
$\Gamma(x)=\int_0^{\infty} e^{-u} u^{x-1} du$
(Theorem \ref{Theorem:Bessel}).
This process is reversible with respect to 
the determinantal point process $\mu_{J_\nu}$ with the {Bessel kernel}
\begin{eqnarray}
K_{J_\nu}(x,y) =
\left\{
   \begin{array}{ll}
\displaystyle{
\frac{J_{\nu}(2 \sqrt{x}) \sqrt{y} J_{\nu}'(2 \sqrt{y})
-\sqrt{x} J_{\nu}'(2\sqrt{x}) J_{\nu}(2\sqrt{y})}{x-y}
} 
& \mbox{if} \quad x \not= y, \\
(J_{\nu}(2 \sqrt{x})^2-J_{\nu+1}(2 \sqrt{x}) J_{\nu-1}(2 \sqrt{x})
& \mbox{if} \quad x=y,
\end{array} \right.
\label{kernel:Bessel}
\end{eqnarray}
with $J_{\nu}'(x)=dJ_{\nu}(x)/dx$.

Finally we will study the 
infinite-dimensional determinantal process 
$(\{\Xi(t)\}_{ t\in [0,\infty)},{\bf P}_{\Ai})$ associated with 
the {\it extended Airy kernel} ${\bK}_{\Ai}$ :
\begin{equation}
{\bK}_{\Ai}(s,x;t,y) \equiv \left\{ 
\begin{array}{cc}
\displaystyle{
\int_{0}^{\infty} d u \, e^{-u(t-s)/2} \Ai(u+x) \Ai(u+y)}
& \mbox{if $t \geq s$} \cr
\displaystyle{- \int_{-\infty}^{0} d u \, e^{-u(t-s)/2} \Ai(u+x) \Ai(u+y)}
& \mbox{if $t < s$},
\end{array}
\right.
\label{def:ex_Airy-kernel}
\end{equation}
$x, y \in\R$, where $\Ai(\cdot)$ is the Airy function defined by
\begin{equation}
\Ai(z) = \frac{1}{2 \pi} \int_{\R} dk \,
e^{i(z k+k^3/3)},
\quad z\in\C,
\nonumber
\end{equation}
which is reversible with respect to 
the determinantal point process $\mu_{\rm Ai}$ with
the Airy kernel
\begin{eqnarray}
\label{kernel:Airy}
K_{\Ai}(x,y) &=&  \left\{ \begin{array}{ll}
\displaystyle{
\frac{\Ai(x) \Ai'(y)-\Ai'(x) \Ai(y)}{x-y}}
&\mbox{if } x \not= y \cr
(\Ai'(x))^2-x (\Ai(x))^2
&\mbox{if } x=y, \end{array} \right.
\end{eqnarray}
where $\Ai'(x)=d \Ai(x)/dx$.
In this process the density $\rho_{\Ai}(x)$ of particles is not bounded
and shows the same asymptotic behavior 
with the function 
$\widehat{\rho}(x) \equiv (\sqrt{-x}/\pi)\1(x<0)$ as $x \to -\infty$,
where ${\bf 1}(\omega)$ is the indicator function of
a condition $\omega$;
${\bf 1}(\omega)=1$ if $\omega$ is satisfied,
and ${\bf 1}(\omega)=0$ otherwise
(see (\ref{asym:Airy}) in Lemma \ref{lemma:correlation_2}).
Therefore the repulsive force from infinite number of particles 
in the negative region $x<0$ causes a positive drift 
with infinite strength. To compensate it 
in order to obtain a stationary system, another negative drift 
with infinite strength should be included in the process. 
For this reason we have to modify our strategy
and we introduce not only a sequence of probability measures
but also a sequence of approximate processes
with finite numbers of particles 
$(\{\Xi_{\widehat{\rho}^N}(t)\}_{t\in [0,\infty)}, \P^{\xi^N}),
N \in \N$, 
having negative drifts such that their strength 
diverges as $N$ goes to infinity.
To determine the $N$-dependence of negative drifts , 
we use the estimate (\ref{asym:AN}) in Lemma \ref{lemma:correlation_2}. 
Then we will show that the limiting process is a version of 
$(\{\Xi(t)\}_{ t\in [0,\infty)},{\bf P}_{\Ai})$ 
and it is Markovian (Theorem \ref{Theorem:Airy}).

\vskip 0.3cm
\noindent{\bf Remark 1.} \quad
Other processes having infinite-dimensional determinantal
point processes as their stationary measures
have been constructed and Markovianity of systems
was proved in Borodin and Olshanski \cite{BO06,BO09,BO10},
Olshanski \cite{Ols10}, and Borodin and Gorin \cite{BG11}.
They determined the state spaces and transition functions
associated with Feller semi-groups and concluded
that these infinite-dimensional determinantal processes
are strong Feller processes.
In the present paper, we first prove Markovianity
of the determinantal processes with the
extended sine kernel $\bK_{\sin}$,
the extended Airy kernel $\bK_{\Ai}$,
and the extended Bessel kernel $\bK_{J_{\nu}}$,
which have been well-studied in the random matrix
theory and its related fields.
To these processes, the argument by \cite{BO06,BO09,BO10,Ols10,BG11}
can not be applied, and the strong Feller property
has not yet proved.
(See also Remark 2 given at the end of Section 2.)
\vskip 0.3cm

The paper is organized as follows. 
In Section 2 preliminaries and main results are given. 
In Section 3 the basic properties of the correlation 
functions are summarized. Section 4 is devoted to proofs of results
by using Lemma 4.3.
In Section 5 Lemma \ref{lemma:correlation_2} is proved 
by using the estimate of Hermite polynomials 
given by Plancherel and Rotach \cite{PR29}.

\SSC{Preliminaries and Main Results}

\subsection{Non-equilibrium dynamics}

For $\xi^N \in \mM$ with $\xi^N(\R)=N\in\N$ and $p 
\in \N_0 \equiv \N \cup \{0\}$
we consider the product
$$
\Pi_{p}(\xi^N, w)= 
\prod_{x \in \supp \xi^{N}}
G \left( \frac{w}{x}, p \right)^{\xi(\{x\})},
\quad w \in \C,
$$
where 
\begin{equation}
G(u,p) = \left\{
   \begin{array}{ll}
\displaystyle{1-u},
& \mbox{if} \quad p=0  \\
& \\
\displaystyle{(1-u)\exp\left [ u+\frac{u^2}{2}+\cdots +\frac{u^p}{p} \right]},
& \mbox{if} \quad p\in\N.
   \end{array} \right. 
\nonumber
\end{equation}
The functions $G(u,p)$ are called
the Weierstrass primary factors \cite{Lev96}.
Then we set
\begin{equation}
\Phi_p(\xi^N, z, w) 
\equiv \Pi_p(\tau_{-z}\xi^N \cap \{0\}^{\rm c}, w-z)
=\prod_{x\in\supp \xi^N \cap \{z\}^{\rm c}}
G\left(\frac{w-z}{x-z},p \right)^{\xi(\{x\})},
\quad w, z \in \C,
\nonumber
\end{equation}
where $\tau_z \xi(\cdot) \equiv \sum_{j \in \Lambda} \delta_{x_j+z}(\cdot)$ 
for $z \in \C$ and $\xi(\cdot) =\sum_{j \in \Lambda} \delta_{x_j}(\cdot)$.


It was proved in Proposition 2.1 of \cite{KT10a} that
the Dyson model $(\{\Xi(t)\}_{ t\in [0,\infty)}, \P^{\xi^N})$,
starting from any fixed configuration $\xi^N\in \mM$ 
is determinantal with the correlation kernel $\mbK^{\xi^N}$
given by 
\begin{eqnarray}
\mbK^{\xi^N}(s, x; t, y)
&=& \frac{1}{2 \pi i} \int_{\R} du \,
\oint_{{\rm C}_{iu}(\xi^N)} dz \,  
p_{\sin}(s, x|z) 
\frac{\Pi_0(\tau_{-z}\xi^N, iu-z)}{iu-z}
p_{\sin}(-t, iu|y)
\nonumber\\
&&  - {\bf 1}(s > t)p_{\sin}(s-t, x|y),
\label{def:K_xi}
\end{eqnarray}
where ${\rm C}_w(\xi^N)$ denotes a closed contour on the
complex plane $\C$ encircling the points in 
$\supp \xi^N$ on the real line $\R$
once in the positive direction but not the point $w$, 
and $p_{\sin}(t, x|y)$ is the generalized heat kernel:
\begin{equation}
p_{\sin}(t, x|y) 
= \frac{1}{\sqrt{2\pi |t|}}\exp\Big\{-\frac{(x-y)^2}{2t}\Big\}{\bf 1}(t\not= 0)
+ \delta(y-x){\bf 1}(t=0),
\quad t \in\R, \ x,y\in\C.
\nonumber
\end{equation}
We put $\mM_0 = \{ \xi \in \mM : \xi(x)\le 1 \, \mbox{for any } x\in\R\}$.
Any element $\xi\in\mM_0$ has no multiple points 
and can be identified with its support which is a countable subset of $\R$.
In case $\xi^N \in \mM_0$, (\ref{def:K_xi}) is rewritten as 
\begin{eqnarray}
\mbK^{\xi^N}(s, x; t, y)
&=& \int_{\R} \xi^N(dx') \, 
\int_{\R} d u \, p_{\sin}(s, x|x')
\Phi_0 (\xi^N, x',iu)p_{\sin}(-t, iu|y)
\nonumber\\
&& - {\bf 1}(s > t) p_{\sin}(s-t, x|y).
\nonumber
\end{eqnarray}

For $\xi \in \mM$, 
$p\in\N_0$, $w, z\in\C$ we define
$$
\displaystyle{
\Phi_p (\xi,z,w)=\lim_{L\to\infty}
\Phi_p(\xi \cap [-L, L], z, w)
}
$$
if the limit finitely exists.
For $L>0, \alpha>0$ and $\xi\in\mM$ we put
$$
M(\xi, L)=\int_{[-L,L]\setminus\{0\}} \frac{\xi(dx)}{x},
\quad \mbox{and} \quad
M(\xi) = \lim_{L\to\infty}M(\xi, L),
$$
if the limit finitely exists,
and for $\alpha > 0$ we put
$$
M_\alpha(\xi)
=\left( \int_{\{0\}^{\rm c}} \frac{1}{|x|^\alpha}\xi(dx)\right)^{1/\alpha}.
$$
It is readily to see that
$\Phi_p (\xi,z,w)$ finitely exists and is not identically $0$
if $M_{p+1}(\xi)<\infty$,
and that $|\Phi_0 (\xi,z,w)| < \infty$ 
if $|M(\xi)|<\infty$ and $M_2( \xi) < \infty$.

For $\kappa>0$, we put
$g^\kappa(x) ={\rm sgn}(x) |x|^{\kappa}$, $x\in\R$,
and $\eta^{\kappa}(\cdot)=\sum_{\ell\in\Z} 
\delta_{g^\kappa(\ell)}(\cdot)$.
For $\kappa\in (1/2,1)$ and $m\in\N$
we denote by $\mY_{\kappa,m}$ the set of configurations $\xi$
satisfying the following conditions 
(C.I) and (C.II):

\vskip 3mm

\noindent {\bf (C.I)}
\qquad\qquad\qquad\qquad\qquad $|M(\xi)|  < \infty$,

\noindent {\bf (C.II)}
\qquad\qquad\qquad\qquad
$\displaystyle{
m(\xi,\kappa)\equiv 
\max_{k\in\Z} \xi\bigg( [ g^\kappa(k), g^\kappa(k+1)] \bigg) \le m.}
$

\noindent And we put
$$
\mY = \bigcup_{\kappa\in (1/2,1)}\bigcup_{m\in\N}\mY_{\kappa,m}.
$$
Noting that the set $\{\xi \in \mM: m(\xi,\kappa)\le m\}$ is 
relatively compact with the vague topology
for each $\kappa\in (1/2,1)$ and $m\in \N$,
we see that $\mY$ is locally compact.
We introduce the following topology on $\mY$.

\begin{df}
\label{def:df2}
Suppose that $\xi, \xi_n \in \mY, n \in \N$.
We say that $\xi_n$ converges $\Phi_0$-moderately to $\xi$, if
\begin{equation}
\lim_{n\to\infty} \Phi_0(\xi_n, i, \cdot)
= \Phi_0(\xi, i, \cdot)
\mbox{ uniformly on any compact set of $\C$.}
\label{eqn:THC2}
\end{equation}
\end{df}

It is easy to see that (\ref{eqn:THC2}) is satisfied, if
$\xi_n$ converges to $\xi$ vaguely
and the following two conditions hold:
\begin{eqnarray}
&&\lim_{L\to\infty} \sup_{n>0 } 
\Bigg| M(\xi_n)- M(\xi_n,L) \Bigg|=0,
\label{eqn:THC3}
\\
&&\lim_{L\to\infty} \sup_{n>0 } 
\Bigg| M_2(\xi_n\cap [-L,L]^c)\Bigg|=0.
\label{eqn:THC4}
\end{eqnarray}
Note that for any $a\in\R$ and $z\in \C$ 
$\displaystyle{\lim_{n\to\infty}\Phi_0(\xi_n,a,z)= \Phi_0(\xi,a,z)}$, 
if $\xi_n$ converges $\Phi_0$-moderately to $\xi$
and $a\notin \supp \xi$.
We denote the space of $\mM$-valued continuous functions defined on $[0,\infty)$ by $C([0,\infty)\to\mM)$. 
We have obtained the following results.
(See Theorem 2.4 of \cite{KT10a} and Theorem 1.4 of \cite{KT10b}.)

\begin{prop}
\label{Theorem:Dyson_CMP}
{\rm (i)} If $\xi\in\mY$, the process 
$(\{\Xi(t)\}_{ t\in [0,\infty)}, \P^{\xi \cap [-L,L]})$
converges to the determinantal process with a correlation kernel $\mbK^{\xi}$
as $L\to\infty$, weakly on the path space $C([0,\infty)\to\mM)$.
In particular, when $\xi \in \mY_{0} \equiv \mY \cap \mM_0$,
$\mbK^{\xi}$ is given by 
\begin{eqnarray}
\mbK^{\xi}(s, x; t, y)
&=& \int_{\R} \xi(dx') \, 
\int_{\R} du \, p_{\sin}(s, x|x')
\Phi_0 (\xi, x',iu)p_{\sin}(-t, iu|y)
\nonumber\\
&& - {\bf 1}(s > t) p_{\sin}(s-t, x|y).
\nonumber
\end{eqnarray}

\noindent {\rm (ii)} 
Suppose that $\xi, \xi_n \in \mY_m^\kappa, n\in\N$,
for some $\kappa\in (1/2,1)$ and $m\in\N$.
If $\xi_n$ converges $\Phi_0$-moderately to $\xi$, then
the process $(\{\Xi(t)\}_{ t\in [0,\infty)}, \P^{\xi_n})$
converges to the process $(\{\Xi(t)\}_{ t\in [0,\infty)}, \P^{\xi})$
as $n\to\infty$, weakly on the path space $C([0,\infty)\to\mM)$.
\end{prop}

\vskip 3mm


Let $\mM^+=\{\xi \cap \R_+ : \xi \in \mM\}$.
We consider a one-parameter family of $\mM^+$-valued processes 
with a parameter $\nu > -1$, 
\begin{equation}
\Xi^{(\nu)}(t, \cdot)=\sum_{j \in \Lambda} \delta_{X^{(\nu)}_j(t)}(\cdot),
\quad t \in [0, \infty),
\label{eqn:Xi1}
\end{equation}
where $X^{(\nu)}_j(t)$'s satisfy the SDEs,
\begin{eqnarray}
&&d X^{(\nu)}_j(t) = 2 \sqrt{X^{(\nu)}_j(t)} dB_j(t)
+2 (\nu+1) dt 
+ 4 X^{(\nu)}_j(t) \sum_{k: k \not=j}
\frac{1}{X^{(\nu)}_j(t)-X^{(\nu)}_k(t)} dt, \quad
\nonumber
\\
&&\qquad\qquad\qquad\qquad\qquad\qquad
j \in \Lambda, \quad t \in [0, \infty)
\label{eqn:ncBESQ}
\end{eqnarray}
and if $-1 < \nu < 0$ with a reflection wall at the origin. 
For a given configuration $\xi^N \in \mM^{+}$ of finite particles, 
$\xi^N(\R_+)=N \in \N$, 
the process starting from $\xi^N$ is denoted by 
$(\{\Xi^{(\nu)}(t)\}_{ t\in [0,\infty)},\P^{\xi^N})$  
and called {\it the noncolliding squared Bessel processes with index $\nu$}.
In Theorem 2.1 of \cite{KT11} it was proved that 
$(\{\Xi^{(\nu)}(t)\}_{ t\in [0,\infty)}, \P^{\xi^N})$ 
is determinantal with the correlation kernel 
\begin{eqnarray}
\mbK^{\xi_N}_{\nu}(s, x; t, y)
&=& \frac{1}{2 \pi i} 
\int_{-\infty}^{0} du
\oint_{{\rm C}_u(\xi^N)} dz  \,  
p^{(\nu)}(s, x|z) \frac{\Pi_0(\tau_{-z}\xi^N, u-z)}{u-z}
p^{(\nu)}(-t, u|y)
\nonumber\\
&& - {\bf 1}(s > t)p^{(\nu)}(s-t, x|y),
\label{eqn:KN1a}
\end{eqnarray}
where for $t \in \R$ and $x, y \in \C$
\begin{eqnarray}
p^{(\nu)}(t, y|x)
&=& \frac{1}{2|t|} \left(\frac{y}{x} \right)^{\nu/2}
\exp \left(-\frac{x+y}{2t} \right)
I_{\nu}\left(\frac{\sqrt{xy}}{|t|} \right)
{\bf 1}(t \not=0, x\not=0)
\nonumber\\
&+&\frac{y^\nu}{(2|t|)^{\nu+1}\Gamma(\nu+1)}
\exp \left(-\frac{y}{2t} \right){\bf 1}(t \not=0, x=0)
+\delta(y-x) {\bf 1}(t=0).
\label{eqn:pnu-}
\end{eqnarray}
Note that for $t\ge 0$ and $x,y\in\R_+$, $p^{(\nu)}(t, y|x)$ is the transition density function of $2(\nu+1)$-dimensional squared Bessel process.
In case $\xi^N \in \mM^+\cap \mM_0$, (\ref{eqn:KN1a}) is equal to
\begin{eqnarray}
\mbK^{\xi^N}_{{\nu}}(s, x; t, y)
&=&  
\int_{0}^{\infty} \xi^N(dx') 
\int_{-\infty}^{0} du \,
 p^{(\nu)}(s, x|x') \Phi_0(\xi^N, x', u)  p^{(\nu)}(-t, u|y)
\nonumber\\
&& - {\bf 1}(s > t)p^{(\nu)}(s-t, x|y)
\nonumber
\end{eqnarray}

For $\kappa \in (1,2)$ and $m\in\N$, let $\mY^+_{\kappa, m}$ be
the set of configurations $\xi$ of $\mM^+$ satisfying 
\noindent {(\bf C.II)},
and put
$$
\mY^+ = \bigcup_{\kappa\in (1,2)}
\bigcup_{m\in\N} \mY^+_{\kappa, m}.
$$
Note that $\mY^+$ is locally compact in the vague topology.
The following proposition is a slight modification of the results
stated in Section 2 of \cite{KT11}, 
which can be proved by the same procedure given there.

\begin{prop}
\label{Theorem:JSP11}
{\rm (i)} If $\xi\in\mY^{+}$, the process
$(\{\Xi^{(\nu)}(t)\}_{ t\in [0,\infty)}, \P^{\xi \cap [-L,L]} )$
converges to the determinantal process 
with a correlation kernel $\mbK^{\xi}_{\nu}$
as $L\to\infty$ in the sense of finite dimensional distributions.
In particular, when 
$\xi \in \mY_{0}^{+} \equiv \mY^{+} \cap \mM_0$,
$\mbK^{\xi}_{\nu}$ is given by 
\begin{eqnarray}
\mbK^{\xi}_{{\nu}}(s, x; t, y)
&=&  
\int_{0}^{\infty} \xi(dx') 
\int_{-\infty}^{0} du \,
 p^{(\nu)}(s, x|x') \Phi_0(\xi_N, x', u)  p^{(\nu)}(-t, u|y)
\nonumber\\
&& - {\bf 1}(s > t)p^{(\nu)}(s-t, x|y)
\nonumber
\end{eqnarray}

\noindent {\rm (ii)} 
Suppose that $\xi, \xi_n \in \mY_{\kappa,m}^{+}, n\in\N$,
for some $\kappa\in (1,2)$ and $m\in\N$.
If $\xi_n$ converges $\Phi_0$-moderately to $\xi$, then
the process $(\{\Xi^{(\nu)}(t)\}_{ t\in [0,\infty)},\P^{\xi_n})$
converges to 
the process $(\{\Xi^{(\nu)}(t)\}_{ t\in [0,\infty)},\P^{\xi})$ as $n \to \infty$ in the sense of finite dimensional distributions.
\end{prop}

\vskip 3mm


Let $\widehat{\rho}^N, N\in\N$ be a sequence of nonnegative functions on $\R$ 
such that $\widehat{\rho}^N(x)=0$ for $x\ge 0$, 
$\int_{\R} dx \ \widehat{\rho}^N(x) =N$ and
$\widehat{\rho}^N(x) \nearrow 
\widehat{\rho}(x) = (\sqrt{-x}/\pi) {\bf 1}(x<0), N \to \infty$.
For $\xi^N\in\mM$ with $\xi^N(\R)=N$ we put
\begin{equation}
M_{\widehat{\rho}^N}(\xi^N)
=\int_{\{0\}^{\rm c}} \frac{\widehat{\rho}^N(x)dx-\xi^N(dx)}{x}.
\nonumber
\end{equation}
and define
\begin{eqnarray}
\Phi_{\widehat{\rho}^N}(\xi^N,w) 
&\equiv& \exp \Bigg[w M_{\widehat{\rho}^N}(\xi^N)\Bigg]
\Pi_1(\xi^N\cap \{0\}^{\rm c}, w), \quad z \in \C, 
\nonumber
\\
\Phi_{\widehat{\rho}^N}(\xi^N,z, w) &\equiv& \Phi_{\widehat{\rho}^N}(\tau_{-z}\xi^N,w-z),
\quad w, z \in \C.
\nonumber
\end{eqnarray}
For $\xi\in\mM$ with $\xi(\R)=\infty$ we put
$\displaystyle{
M_{\cA}(\xi)=\lim_{N\to\infty}M_{\widehat{\rho}^N}(\xi^N),
}$
if $\xi^N$ converges to $\xi$, $N\to \infty$, with the vague topology
and
\begin{equation}
\lim_{L\to\infty}\sup_{N\in\N}\int_{|x|>L} 
\frac{\widehat{\rho}^N(x)dx-\xi^N(dx)}{x}=0,
\label{condition:convergence}
\end{equation}
is satisfied. Then $M_{\cA}(\xi)$ is also represented as
\begin{equation}
M_{\cA}(\xi)
=\lim_{L\to\infty} \int_{0<|x|<L} \frac{\widehat{\rho}(x)dx-\xi(dx)}{x}.
\label{def:MA}
\end{equation}
For $\xi \in \mM$ with $M_{\cA}(\xi)<\infty$ and $z\in\C$, we define
\begin{eqnarray}
&&\Phi_{\cA}(\xi,w) 
\equiv \exp \Bigg[ w M_{\cA}(\xi) \Bigg]
\Pi_1(\xi\cap \{0\}^{\rm c}, w), \quad w \in \C, 
\nonumber\\
&&\Phi_{\cA} (\xi,z,w)\equiv \Phi_{\cA}(\tau_{-z}\xi,w-z), \quad w, z \in \C.
\nonumber
\end{eqnarray}
We note that $\Phi_{\cA} (\xi,z,w)$ finitely exists and $\Phi_{\cA}(\xi,z,w) \not \equiv0$,
if $|M_{\cA}(\xi)|<\infty$ and $M_2(\xi) < \infty$.

Consider the process with a finite number of particles given by 
$\Xi_{\widehat{\rho}^N}(t)=\sum_{j=1}^{N} \delta_{Y_j(t)}$ with
\begin{equation}
Y_j(t) = X_j(t)+\frac{t^2}{4} +t \int_{\R}\frac{\widehat{\rho}^N(x)dx}{x}, \quad
1 \leq j \leq N, \quad
t \in [0, \infty),
\label{eqn:Y_j}
\end{equation}
associated with the solution
$\X(t)=(X_1(t), \dots, X_N(t))$ of (\ref{eqn:Dyson}).
In other words, $\Y(t)=(Y_1(t), Y_2(t), \dots, Y_N(t))$
satisfies the following SDEs ;
\begin{eqnarray}
&&dY_j(t) = dB_j(t)+\left( \frac{t}{2}+\int_{\R}\frac{\widehat{\rho}^N(x)dx}{x} \right) dt
+\sum_{\substack{1 \leq k \leq N \\ k \not= j}}\frac{dt}{Y_j(t)-Y_k(t)},
\nonumber\\
&&\qquad\qquad\qquad\qquad\qquad\qquad\qquad\qquad
\ 1 \leq j \leq N, \ t \in [0, \infty).
\nonumber
\end{eqnarray}
It is proved in Proposition 2.4 of \cite{KT09} that the process 
$(\{\Xi_{\widehat{\rho}^N}(t)\}_{t\in [0,\infty)},\P^{\xi^N})$ starting from any fixed configuration $\xi^N$ with a finite number of particles is determinantal with the correlation kernel $\mbK_{\widehat{\rho}^N}^{\xi^N}$
given by 
\begin{eqnarray}
\mbK^{\xi^{N}}_{\widehat{\rho}^N}(s,x ;t,y)
&=& \frac{1}{2\pi i}
\int_{\R} du \,
\oint_{{\rm C}_{iu}(\xi_N)} dz \,
q(0,s, x-z)
\frac{\Pi_0(\tau_{-z}\xi,iu-z)}{iu-z}
\nonumber\\
&&\qquad\times
\exp\bigg[(iu-z)\int_{\R}\frac{\widehat{\rho}^N(v)dv}{v}\bigg]
q(t,0, iu-y)
\nonumber\\
&& - {\bf 1}(s>t)q(t, s, x-y),
\label{eqn:K3}
\end{eqnarray}
where $q(s,t,y-x), s, t \in \R, s \not= t, x, y \in \C$ is given by
\begin{eqnarray}
&& q(s, t, y-x) =
p_{\sin}\left(t-s, \left(y-\frac{t^2}{4} \right)-
\left( x-\frac{s^2}{4} \right) \right)
\nonumber\\
\label{eqn:q}
&& = \frac{1}{\sqrt{2 \pi |t-s|}}
\exp \left[ 
-\frac{(y-x)^2}{2(t-s)}+\frac{(t+s)(y-x)}{4}
-\frac{(t-s)(t+s)^2}{32} \right].
\nonumber
\end{eqnarray}
Note that $q(s,t,y-x)$, $t>s\ge 0$, $x,y \in\R$ 
is the transition density function of the process 
$B(t)+t^2/4$,
where $B(t), t \in [0, \infty)$ is the one-dimensional
standard Brownian motion.

In case $\xi^N \in \mM_0$, (\ref{eqn:K3}) is equal to
\begin{eqnarray}
&& \mbK^{\xi^N}_{\widehat{\rho}^N}(s,x;t,y)
=\int_{\R}\xi(dx') 
\int_{\R} du \, 
q(0,s, x-x') 
\Phi_{\widehat{\rho}^N}(\xi^N, x', iu)
q(t, 0, iu-y)
\nonumber\\
&& \qquad\qquad\qquad - {\bf 1}(s>t)q(t, s, x-y).
\nonumber
\end{eqnarray}

For $\kappa\in (1/2,2/3)$ and $m\in\N$, we denote by 
$\mY^{\cA}_{\kappa, m}$ 
the set of configurations $\xi$
satisfying the conditions

\vskip 3mm

\noindent ({\bf C.I-$\cA$})
\qquad\qquad\qquad\qquad 
$|M_{\cA}(\xi)|  < \infty$,

\vskip 3mm

\noindent 
and ({\bf C.II}). And we define the space of configurations 
$$
\mY^{\cA} = \bigcup_{\kappa\in (1/2,2/3)}\bigcup_{m\in\N}
\mY^{\cA}_{\kappa, m}.
$$
Note that $\mY^{\cA}$ is locally compact.
We introduce the following topology on $\mY^{\cA}$.

\begin{df}
\label{def:dfPA_con}
Suppose that $\xi, \xi_n \in \mY^{\cA}, n \in \N$.
We say that $\xi_n$ converges $\Phi_\cA$-moderately to $\xi$, if
\begin{equation}
\lim_{n\to\infty} \Phi_{\cA}(\xi_n, i, \cdot)
= \Phi_{\cA}
(\xi, i, \cdot)
\mbox{ uniformly on any compact set of $\C$.}
\label{eqn:PA}
\end{equation}
\end{df}


The following proposition is a slight modification of the results
stated in Section 2.4 of \cite{KT09}. 

\begin{prop}
\label{Theorem:JSP}
\quad Let $\xi\in\mY^{\cA}$ and $\xi^N, N\in\N$ be a sequence of configurations such that $\xi^N(\R)=N$ and $\xi^N$ converges to $\xi$ with the vague topology. Suppose that {\rm (\ref{condition:convergence})} is satisfied and 
$\displaystyle{\max_{N\in\N}m(\xi^N,\kappa)\le m}$,
for some $\kappa\in (1/2,2/3)$, $m\in\N$.

\noindent {\rm (i)} \,
The family of the distributions of the processes 
$(\{\Xi_{\widehat{\rho}^N}(t)\}_{t\in [0,\infty)}, \P^{\xi^N})$
is tight in the space of probability measures on $C([0,\infty)\to\mM)$.

\noindent {\rm (ii)} \,
The sequence of the processes $(\{\Xi_{\widehat{\rho}^N}(t)\}_{t\in [0,\infty)}, \P^{\xi^N})$
converges to the determinantal process $(\{\Xi_{\cA}(t)\}_{t\in [0,\infty)}, \P^{\xi})$ with a correlation kernel $\mbK^{\xi}_{\cA}$ as $N\to\infty$, 
weakly on the path space $C([0,\infty)\to\mM)$.
In particular, when 
$\xi \in \mY_{0}^{\cA} \equiv \mY^{\cA} \cap \mM_0$,
$\mbK^{\xi}_{\cA}$ is given by 
\begin{eqnarray}
&& \mbK^{\xi}_{\cA}(s,x;t,y)
=\int_{\R}\xi(dx') 
\int_{\R} du \, 
q(0,s, x-x') 
\Phi_{\cA}(\xi, x', iu)
q(t, 0, iu-y)
\nonumber\\
&& \qquad\qquad\qquad - {\bf 1}(s>t)q(t, s, x-y).
\nonumber
\end{eqnarray}

\noindent {\rm (iii)} 
Suppose that $\xi, \xi_n \in \mY_{\kappa,m}^{\cA}, n\in\N$,
for some $\kappa\in (1/2,2/3)$ and $m\in\N$.
If $\xi_n$ converges $\Phi_{\cA}$-moderately to $\xi$, then
the process $(\{\Xi_{\cA}(t)\}_{ t\in [0,\infty)},\P^{\xi_n})$
converges to 
the process $(\{\Xi_{\cA}(t)\}_{ t\in [0,\infty)},\P^{\xi})$ as $n \to \infty$weakly on the path space $C([0,\infty)\to\mM)$.
\end{prop}

By the same argument in \cite{KT10b} we can obtain the first assertion (i), under which the other assertions (ii) and (iii) can be proved by the same procedure given in Section 2.4 of \cite{KT09}. Hence, we skip the proof of this proposition.

\subsection{Main results}

For $\xi\in \mY$ we put
\begin{equation}
T_t f(\xi) =\E^{\xi} \Big[f(\Xi(t)) \Big], \quad t\ge 0,
\label{eqn:defTt}
\end{equation}
for a bounded vaguely continuous function $f$ on $\mM$,
where $\E^{\xi}$ represents the expectation with respect to the probability measure $\P^{\xi}$.
When $\xi, \xi_n \in \mY_{\kappa, m}, n\in\N$,
for some $\kappa\in (1/2,1)$ and $m\in\N$,
$T_t f(\xi_n)$ converges to $T_t f(\xi)$,
if $\xi_n$ converges $\Phi_0$-moderately to $\xi$, as $n \to \infty$.
We denote by $L^2(\mM,{\mu})$ 
the space of square integrable functions on $\mM$
with respect to the probability measure ${\mu}$,
which is equipped with the inner product
$\displaystyle{\langle f, g \rangle_{\mu}
\equiv \int_{\mM} \mu(d \xi) f(\xi) g(\xi),
f,g \in L^2(\mM,{\mu})}$.
We write the expectation with respect to the probability measure ${\bf P}_{\sin}$ as ${\bf E}_{\sin}$.
The first main theorem of the present paper is the following.

\begin{thm}
\label{Theorem:Dyson}
\noindent {\rm (i)} 
${\mu_{\sin}} (\mY)=1$ and 
$T_t$ is extended to the contraction 
operator on $L^2(\mM,{\mu_{\sin}})$.

\noindent {\rm (ii)} 
The processes $(\{\Xi(t)\}_{ t \in [0, \infty)}, \P_{\mu_{\sin}})$
and $(\{\Xi(t)\}_{t \in [0, \infty)},{\bf P}_{\sin})$ 
are versions of the same determinantal process.
In particular, for any $t\ge 0$
\begin{equation}
{\bf E}_{\sin} \Big[f_0(\Xi(0))f_1(\Xi(t))\Big]
= \langle f_0, T_t f_1 \rangle_{\mu_{\sin}},
\quad f_0, f_1 \in L^2(\mM,{\mu_{\sin}}).
\label{eqn:inner_product}
\end{equation}

\noindent {\rm (iii)} 
The reversible process $(\{\Xi(t)\}_{t\in [0,\infty)}, {\bf P}_{\sin})$
is Markovian.
\end{thm}

\vskip 3mm

For $\xi\in \mY^{+}$ we put
\begin{equation}
T^{(\nu)}_t f(\xi) =\E^{\xi} \Big[f(\Xi^{(\nu)}(t)) \Big],
 \quad t\ge 0,
\label{eqn:defTtB}
\end{equation}
for a bounded continuous function $f$ on $\mM^{+}$.
When $\xi, \xi_n \in \mY^{+}_{\kappa,m}, n\in\N$,
for some $\kappa\in (1,2)$ and $m\in\N$,
$T^{(\nu)}_t f(\xi_n)$ converges to $T^{(\nu)}_t f(\xi)$,
if $\xi_n$ converges $\Phi_0$-moderately 
to $\xi$, as $n \to \infty$.
We write the expectation with respect to the probability measure ${\bf P}_{J_\nu}$ as ${\bf E}_{J_\nu}$.
The second main theorem of the present paper is the following.

\begin{thm}
\label{Theorem:Bessel}
\noindent {\rm (i)} 
${\mu_{J_\nu}} (\mY^{+})=1$ and 
$T^{(\nu)}_t$ is extended to the contraction 
operator on $L^2(\mM^+,{\mu_{J_\nu}})$.

\noindent {\rm (ii)} 
The processes $(\{\Xi^{(\nu)}(t)\}_{t \in [0, \infty)}, \P_{\mu_{J_\nu}})$
and $(\{\Xi(t)\}_{t \in [0, \infty)}, {\bf P}_{J_\nu})$ 
are versions of the same determinantal process.
In particular, for any $t\ge 0$
\begin{equation}
{\bf E}_{J_\nu} \Big[f_0(\Xi(0))f_1(\Xi(t))\Big]
= \langle f_0, T^{(\nu)}_t f_1 \rangle_{\mu_{J_\nu}},
\quad f_0, f_1 \in L^2(\mM^+,{\mu_{J_\nu}}).
\label{eqn:inner_product_B}
\end{equation}

\noindent {\rm (iii)} 
The reversible process $(\{\Xi(t)\}_{t\in [0,\infty)}, {\bf P}_{J_\nu})$
is Markovian.
\end{thm}

\vskip 3mm


For $\xi\in \mY^{\cA}$ we put
\begin{equation}
T^{\cA}_t f(\xi) =\E^{\xi} \Big[f(\Xi_{\cA}(t)) \Big],  
\quad t\ge 0,
\label{eqn:defTtA}
\end{equation}
for a bounded continuous function $f$ on $\mM$.
When $\xi, \xi_n \in \mY^{\cA}_{\kappa,m}, n\in\N$,
for some $\kappa\in (1/2,2/3)$ and $m\in\N$,
$T^{\cA}_t f(\xi_n)$ converges to $T^{\cA}_t f(\xi)$,
if $\xi_n$ converges $\Phi_{\cA}$-moderately 
to $\xi$, as $n \to \infty$.
We write the expectation with respect to the probability measure ${\bf P}_{\Ai}$ as ${\bf E}_{\Ai}$.
The last main theorem of the present paper is the following.

\begin{thm}
\label{Theorem:Airy}
\noindent {\rm (i)} 
${\mu_{\Ai}} (\mY^{\cA})=1$ and 
$T^{\cA}_t$ is extended to the contraction 
operator on $L^2(\mM,{\mu_{\Ai}})$.

\noindent {\rm (ii)} 
The processes $(\{\Xi_{\cA}(t)\}_{ t \in [0, \infty)}, \P_{\mu_{\rm Ai}})$
and
$(\{\Xi(t)\}_{t \in [0, \infty)}, {\bf P}_{\rm Ai})$
are versions of the same determinantal process. 
In particular, for any $t\ge 0$
\begin{equation}
{\bf E}_{\Ai} \Big[f_0(\Xi(0))f_1(\Xi(t))\Big]
= \langle f_0, T^{\cA}_t f_1 \rangle_{\mu_{\Ai}},
\quad f_0, f_1 \in L^2(\mM,{\mu_{\Ai}}).
\label{eqn:inner_product_A}
\end{equation}

\noindent {\rm (iii)} 
The reversible process $(\{\Xi(t)\}_{t\in [0,\infty)}, {\bf P}_{\Ai})$
is Markovian.
\end{thm}

\noindent{\bf Remark 2.} \quad
A function $f$ on the configuration space $\mM$ is said to be
polynomial, if it is written in the form
$f(\xi)=F \left(
\int_{\R}\phi_1(x)\xi(dx), \int_{\R}\phi_2(x)\xi(dx),
\dots,\int_{\R}\phi_k(x)\xi(dx) \right)$
with a polynomial function $F$
on ${\R}^{k}, k \in \N$,
and smooth functions $\phi_j, 1 \leq j \leq k$
on $\R$ with compact supports.
Let $\wp$ be the set of all polynomial functions on $\mM$.
In \cite{KT07b} we introduced a class of determinantal reversible processes including the three processes 
$(\{\Xi(t)\}_{t \in [0, \infty)},{\bf P}_{\sin})$,
$(\{\Xi(t)\}_{t \in [0, \infty)}, {\bf P}_{J_\nu})$ and
$(\{\Xi(t)\}_{ t \in [0, \infty)}, {\bf P}_{\Ai})$,
and showed in Proposition 7.2 that
for any element $(\{{\Xi}(t)\}_{t \in [0, \infty)},{\bf P}_\mu)$
of the class with the reversible measure $\mu$,
the following relation holds:
\begin{eqnarray}
\lim_{t\to 0}\frac{-1}{t}{\bf E}_\mu
\big[ f(\Xi(0))g(\Xi(t))\big]
&=& \int_{\R}\rho(x)dx
\int_{\mM} \mu^x(d\eta)
\frac{\partial}{\partial x}f(\eta+\delta_{x})
\frac{\partial}{\partial x}g(\eta+\delta_{x})
\nonumber\\
&\equiv& {\cal E}_0^\mu(f,g),
\quad \mbox{ for any } f,g \in \wp,
\nonumber
\end{eqnarray}
where $\rho$ is the one-point correlation function of $\mu$ and
$\mu^{x}$ is the {\it Palm measure} of $\mu$, satisfying
the following relation with $\mu$, 
$$
\int_{\mM}\mu(d\xi) \xi(K)f(\xi)
=\int_K dx \rho(x) \int_{\mM}\mu^x (d\eta) f(\eta+\delta_x)
$$
for any $f\in\wp$ and any compact subset $K$ of $\R$.
Theorems \ref{Theorem:Dyson}, \ref{Theorem:Bessel} and \ref{Theorem:Airy}
imply that the bilinear forms 
$({\cal E}_0^{\mu_{\sin}}, \wp)$, $({\cal E}_0^{\mu_{J_\nu}}, \wp)$ 
and $({\cal E}_0^{\mu_{\Ai}}, \wp)$ are pre-Dirichlet forms \cite{FOT94} for the three process 
$(\{\Xi(t)\}_{ t\in [0,\infty)},{\bf P}_{\sin})$,
$(\{\Xi(t)\}_{ t\in [0,\infty)},{\bf P}_{J_\nu})$ and $(\{\Xi(t)\}_{ t\in [0,\infty)}, {\bf P}_{\Ai})$, respectively.

Spohn \cite{Spo87} considered an infinite particle system
obtained by taking the $N \to \infty$ limit of 
(\ref{eqn:Dyson}) and 
studied the equilibrium dynamics with respect
to the determinantal point process $\mu_{\sin}$.
From the Dirichlet form approach
Osada \cite{Osa96} constructed the infinite particle systems 
represented by the diffusion process associated 
with the Dirichlet form which is the minimal closed 
extension of $({\cal E}_0^{\mu_{\sin}}, \wp)$.
Recently he proved that this system satisfies
the SDEs (\ref{eqn:Dyson}) with $N= \infty$ \cite{Osa08}.
The equivalence of the determinantal processes 
$(\{\Xi(t)\}_{ t\in [0,\infty)},{\bf P}_{\sin})$ 
and the infinite-dimensional equilibrium dynamics of 
Spohn and Osada is, however, not yet proved.
It is relevant to the uniqueness of Markov extensions of 
the pre-Dirichlet form associated 
with an infinite particle system with long-range interaction, 
which is an interesting open problem.
For an infinite particle system with finite-range interaction,
there are some results on the uniqueness problem (see, e.g., \cite{T97}).

\SSC{Correlation functions}

\subsection{Some estimates for determinantal point processes}

Let $\mu$ be a probability measure on $\mM$
with correlation functions 
$\rho_m(\x_m)$, $\x_m \in \R^m$,
$m\in\N$.
Then for $f\in\C_0(\R)$ and $\theta\in\R$
\begin{eqnarray}
\Psi(f,\theta)
\equiv \int_{\mM}\mu(d\xi) e^{\theta\int_{\R^d}f(x)\xi(dx)}
=\sum_{m=0}^\infty \frac{1}{m !}
\int_{\R^{m}} \prod_{j=1}^{m} d x_{j}
\prod_{j=1}^m \Big(e^{\theta f(x_j)}-1\Big)
\rho_m \Big( \x_m \Big).
\nonumber
\end{eqnarray}
Let ${\cal K}$ be a symmetric linear operator with kernel $K$.
In this section we assume that the operator ${\cal K}$  satisfies 
{\bf Condition A} in \cite{ST03}.
The probability measure $\mu$ is called a determinantal point process
with correlation kernel $K$,
if its correlation functions are represented by
$$
\rho_m(\x_m)
=\det_{
{1 \leq j,k \leq m}
}
\Big[
K(x_j,x_k)
\Big].
$$
We often write $\rho$ for $\rho_1$, and $\rho(A)$ 
for $\int_A \rho(x) dx$, $A\in{\cal B}(\R)$, respectively.
In this subsection we give some lemmas which will 
be used in Section \ref{sec4}.

\begin{lem}
\label{lemma:correlation_1}
Let $\mu$ be a determinantal point process. Then
for any bounded closed interval $D$ of $\R$, we have
\begin{equation}
\int_{\mM} \mu(d\eta) \Big| \eta(D)-\int_D \rho (x)dx\Big|^{2k}
\le ( 3 \rho(D) )^k, \quad k\in\N.
\label{est:kth}
\end{equation}
\end{lem}

\noindent{\it Proof.} Let $K$ be the correlation kernel
for the determinantal point process $\mu$.
We put $K_D(x,y) = \1_D(x)K(x,y)\1_D(y)$
with $\1_D(x) \equiv \1(x \in D)$ and 
introduce the linear operator $K_D$ defined by
$$
{\cal K}_Df(x) = \int_{\R} K_D(x,y)f(y)dy, \quad f\in L^2(\R,dx).
$$
Since $D$ is compact, $K_D$ is in the trace class. 
Then nonzero spectrums of the operator are 
eigenvalues $\{\lambda_j\}_{j\in\Lambda}$ with $0< \lambda_j \le 1$,
$j\in\Lambda$.
Hence, the distribution of $\eta(D)$ under $\mu$
is identical with the distribution of the random variable
$\sum_{j\in\Lambda}X_j$,
where $X_j$, $j\in\Lambda$ are independent 
Bernoulli random variables with 
$P(X_j=1)=\lambda_j$ and $P(X_j=0)=1-\lambda_j$
(see for instance \cite{ST03}).
In fact we can readily see that
\begin{eqnarray}
\int_{\mM}e^{\theta \eta(D)} \mu(d\eta)
&=& \Det \left[ 
\delta(x-y) + K(x,y)(e^{\theta\1_D(x)}-1)
\right]
\nonumber\\
&=&\Det \left[ 
\delta(x-y) + K_D(x,y)(e^\theta-1)
\right]
= \prod_{j\in\Lambda}\left\{
1+\lambda_j(e^\theta-1) \right\},
\nonumber
\end{eqnarray}
which coincides with the generating function of 
$\sum_{j\in\Lambda}X_j$.
By this identity of distributions we have
\begin{eqnarray}
\widehat{\Psi}(\theta)
&\equiv& \int_{\mM}\mu(d\eta) 
\exp\left\{
\theta \Big( \eta(D)-\int_{\mM}\mu(d\eta)\eta(D)\Big)^2
\right\}
\nonumber\\
&=& E\left[
\exp \left\{\theta \sum_{j\in\Lambda}(X_j -EX_j)^2
\right\}\right]
=\prod_{j\in\Lambda}e^{\theta \lambda^2}
\left\{ 1-\lambda_j (1-e^{-2\theta\lambda_j})\right\}.
\nonumber
\end{eqnarray}
By simple calculations with the relation $0\le \lambda_j\le 1$,
we see that, for any $k\in\N$
$$
\frac{d^k}{d\theta^k}e^{\theta\lambda_j^2}
\le \frac{d^k}{d\theta^k}e^{\theta\lambda_j},
\quad 
\frac{d^k}{d\theta^k}\{ 1-\lambda_j+\lambda_je^{-2\theta\lambda_j}\}
\le \frac{d^k}{d\theta^k}e^{2\theta \lambda_j},
\quad \theta \ge 0.
$$
Hence,
\begin{eqnarray}
&&\int_{\mM} \mu(d\eta) \Big| \eta(D)-\int_{\mM}\mu(d\eta)\eta(D)\Big|^{2k}
= \frac{d^k}{d\theta^k}\widehat{\Psi}(\theta)\Big|_{\theta=0}
\nonumber\\
&&\qquad\qquad\qquad \le \left. \frac{d^k}{d\theta^k}
\exp \left\{3\theta \sum_{j\in\Lambda}\lambda_j \right\}
\right|_{\theta=0}
\le \left(3\sum_{j\in\Lambda}\lambda_j\right)^k.
\nonumber
\end{eqnarray}
Since $\int_D \rho (x)dx=\int_{\mM}\mu(d\xi)\eta(D)
=E\left[\sum_{j\in \Lambda}X_j\right]= \sum_{j\in\Lambda}\lambda_j$,
we obtain the lemma. \qed
\vskip 3mm

\begin{lem}
\label{lemma:correlation_3}
Let $\xi \in \mM$ and $\rho$ be the nonnegative function on $\R$.
Suppose that there exist $\varepsilon \in (0,1)$, 
$C_1>0$ and $L_1>0$ such that
\begin{equation}
\label{eqn:72}
\left|\xi([0,L])-\int_0^L \rho(x)dx \right|\le C_1L^{\varepsilon},
\quad
\left|\xi([-L,0))-\int_{-L}^0 \rho(x)dx \right|\le C_1L^{\varepsilon},
\quad L\ge L_1.
\end{equation}

\noindent {\rm (i)}
If 
\begin{equation}
\label{eqn:71_i}
\lim_{L\to\infty}\int_{1\le |x| \le L} \frac{\rho(x) dx}{x}
\quad \mbox{finitely exits},
\end{equation}
then $\xi$ satisfies {\bf (C.I)}.

\noindent {\rm (ii)}
If a nonnegative function $\overline{\rho}$ on $\R$ satisfies 
\begin{equation}
\label{eqn:71_ii}
\left| \int_{|x|\ge L} 
\frac{\rho(x) dx-\overline{\rho}(x)dx}{x}
\right|
\le C_2 L^{-\delta},
\end{equation}
for some $\delta>0$, then $\xi$ satisfies
\begin{equation}
\bigg| \int_{|x|\ge L} \frac{\overline{\rho}(x)dx- \xi(dx)}{x}\bigg|\le
\frac{3C_1}{1-\varepsilon}L^{\varepsilon -1} +C_2L^{-\delta}.
\end{equation}
\end{lem}

\noindent{\it Proof.} 
By using the integration by parts formula with (\ref{eqn:72}),
we see that for $L_2 > L \ge L_1$
$$
\bigg| \int_{L}^{L_2} \frac{\rho(x) dx}{x} 
- \int_{L}^{L_2} \frac{\xi(dx)}{x} \bigg|
\le 2C_1L^{\varepsilon -1}+C_1 \int_{L}^{L_2} x^{\varepsilon -2}dx
\le \frac{C_1}{1-\varepsilon}
(3L^{\varepsilon-1} - {L_2}^{\varepsilon-1}).
$$
Similarly, we have
$$
\bigg| \int_{-L_2}^{-L} \frac{\rho(x)dx}{x} 
- \int_{-L_2}^{-L} \frac{\xi(dx)}{x} \bigg|
\le \frac{C_1}{1-\varepsilon}
(3 L^{\varepsilon-1} - {L_2}^{\varepsilon-1}).
$$
Then for $L\ge L_1$
$$
\bigg| \int_{|x|\ge L} \frac{\xi(dx)}{x}
- \int_{|x|\ge L} \frac{\rho(x)dx}{x}\bigg|
\le \frac{3C_1}{1-\varepsilon}L^{\varepsilon-1}.
$$
Then (i) and (ii) are derived easily.
\qed

\vskip 3mm

\begin{lem}
\label{lemma:correltion_4}
Let $\mu$ be a probability measure on $\mM$ with 
correlation functions $\rho_m$, $m\in\N$.
Suppose that there exist $m'\in\N$ and $p<m'-1$
such that
\begin{eqnarray}
\label{eqn:74}
&&\int_{\mM} \mu(d\xi) \Big| \xi([0,L))-\int_0^L \rho (x)dx \Big|^{m'}
= {\cal O}(L^p), \quad L\to\infty,
\\
\label{eqn:75}
&&\int_{\mM} \mu(d\xi) \Big| \xi([-L,0))-\int_{-L}^0 \rho (x)dx \Big|^{m'}
= {\cal O}(L^p), \quad L\to\infty.
\end{eqnarray}

\noindent {\rm (i)}
In case {\rm (\ref{eqn:71_i})} holds and
\begin{equation}
\label{eqn:73}
\sum_{k\in\Z}
\int_{[g^{\kappa}(k),g^{\kappa}(k+1)]^m}\rho_m(\x_m)d\x_m  <\infty
\end{equation}
is satisfied for some $m\in\N$ and $\kappa \in (1/2, 1)$, then $\mu (\mY)=1$.

\noindent {\rm (ii)}
In case $\mu (\mM^+)=1$ and {\rm (\ref{eqn:73})} is
satisfied for some $m\in\N$ and $\kappa \in (1, 2)$,
then $\mu (\mY^{+})=1$.

\noindent {\rm (iii)}
In case {\rm (\ref{eqn:71_ii})} holds with 
$\overline{\rho}=\widehat{\rho}$ and {\rm (\ref{eqn:73})} 
satisfied for some $m\in\N$ and $\kappa \in (1/2, 2/3)$,
then $\mu (\mY^{\cA})=1$.
\end{lem}

\noindent{\it Proof.} 
Take $\varepsilon \in ((p+1)/m',1)$.
Using Chebyshev's inequality with (\ref{eqn:74}) and (\ref{eqn:75}),
we can find a positive constant $C$ such that
$$
\mu \left(|\xi((0,L])- \int_0^L \rho (x)dx|\ge CL^{\varepsilon} \right)
\le CL^{p-m'\varepsilon},
$$
and
$$
\mu \left(|\xi([-L,0))- \int_{-L}^0 \rho (x)dx|\ge CL^{\varepsilon} \right)
\le CL^{p-m'\varepsilon}
$$
Since $p-m'\varepsilon <-1$, we have
\begin{eqnarray}
&&\sum_{L\in\N}\mu \left(|\xi((0,L])-\int_0^L \rho (x)dx| )
|\ge CL^{\varepsilon} \right)<\infty,
\nonumber\\
&&\sum_{L\in\N}\mu \left(|\xi([-L,0))- \int_{-L}^0 \rho (x)dx|
\ge CL^{\varepsilon} \right)<\infty,
\nonumber
\end{eqnarray}
and so the condition (\ref{eqn:72}) is satisfied $\mu$-a.s. $\xi$
by the Borel-Cantelli lemma.
Since we have assumed the condition (\ref{eqn:71_i}) in (i), 
we can apply Lemma \ref{lemma:correlation_3} and 
conclude that ({\bf C.I}) holds for $\mu$-a.s. $\xi$.
Similarly, since the condition (\ref{eqn:71_ii}) 
with $\overline{\rho}=\widehat{\rho}$ is provided in (iii), 
we conclude that ({\bf C.I-$\cA$}) holds for $\mu$-a.s. $\xi$.

The condition (\ref{eqn:73}) implies
$$
\sum_{k\in\Z}\mu \Big(\xi(g^{\kappa}(k), g^{\kappa}(k+1))>m \Big) 
<\infty.
$$
Then, by the Borel-Cantelli lemma again,
({\bf C.II}) is derived.
Therefore, we obtain the lemma.
\qed

\subsection{Equilibrium dynamics}

The Dyson model starting from $N$ points all at the origin 
is determinantal with the correlation kernel
\begin{eqnarray}
&&\mbK_N(s,x;t,y) = \left\{
   \begin{array}{ll}
\displaystyle{
\frac{1}{\sqrt{2s}}
\sum_{k=0}^{N-1}
\left(\frac{t}{s}\right)^{k/2}
\varphi_{k}\left(\frac{x}{\sqrt{2s}}\right)
\varphi_{k}\left(\frac{y}{\sqrt{2t}}\right)}
& \quad \mbox{if} \ s\leq t \\
\displaystyle{\frac{-1}{\sqrt{2s}}
\sum_{k=N}^{\infty} 
\left(\frac{t}{s}\right)^{k/2}
\varphi_{k}\left(\frac{x}{\sqrt{2s}}\right)
\varphi_{k}\left(\frac{y}{\sqrt{2t}}\right)}
& \quad \mbox{if} \ s > t,
   \end{array} \right. 
\label{def:K_finite}
\end{eqnarray}
where $h_k=\sqrt{\pi}2^k k!$ and 
\begin{equation}
\varphi_k(x)=\frac{1}{\sqrt{h_k}}e^{-x^2/2}H_k(x),
\label{def:phin}
\end{equation}
with the Hermite polynomials $H_k, k\in \N_0$ 
(see Eynard and Mehta \cite{EM98}).
The distribution of the process at time $t>0$ 
is identical with $\mu^{\rm GUE}_{N,t}$, the distribution of eigenvalues
in the GUE with variance $t$.
The extended sine kernel ${\bf K}_{\sin}$ with density 1 
is derived as the bulk scaling limit \cite{NF98}:
\begin{equation}
\mbK_N \left( \frac{2N}{\pi^2} +s, 
x; \frac{2N}{\pi^2} +t, y \right)
\to {\bf K}_{\sin}(s,x;t,y),
\quad N\to\infty,
\label{limit:bulk}
\end{equation}
which implies
\begin{eqnarray}
&& (\{\Xi(t)\}_{ t \in [0,\infty)}, \P_{\mu^{\rm GUE}_{N,2N/\pi^2}})
\, \fd \,
(\{\Xi(t)\}_{t\in [0,\infty)}, {\bf P}_{\sin}), \quad N \to \infty,
\nonumber\\
\label{conv:vaguely1}
&& \hskip 3cm 
\quad \mbox{in the vague topology},
\end{eqnarray}
where $\fd $ means the convergence in the sense of finite dimensional distributions.

The noncolliding squared Bessel process 
starting from $N$ points all at the origin
is determinantal with the correlation kernel 
$$
\mbK^{(\nu)}_{N}(s,x;t,y)= \left\{
   \begin{array}{ll}
\displaystyle{
\frac{1}{2s}\sum_{k=0}^{N-1}
\left(\frac{t}{s}\right)^k
\varphi_k^\nu \left(\frac{x}{2s}\right)
\varphi_k^\nu \left(\frac{y}{2t}\right),}
& \quad \mbox{if $s \leq t$,}
\\
\displaystyle{
-\frac{1}{2s}\sum_{k=N}^{\infty}
\left(\frac{t}{s}\right)^k
\varphi_k^\nu \left(\frac{x}{2s}\right)
\varphi_k^\nu \left(\frac{y}{2t}\right),}
& \quad \mbox{if $s > t$,}
\end{array} \right. 
$$
where $\varphi_k^\nu (x) = 
\sqrt{\Gamma(k+1)/\Gamma(\nu+k+1)}x^{\nu/2}
L_k^{\nu}(x)e^{-x/2}$ 
with the Laguerre polynomials $L_k^{\nu}(x), k\in\N_0$
with parameter $\nu>-1$ \cite{KT07a}.
Let $\mu^{(\nu)}_{N,t}$ be
the distribution of the process at time $t>0$.
In the case that $\nu\in \N_0$, the probability measure 
$\mu^{(\nu)}_{N,t}$ is identical with the distribution of 
eigenvalues in the chGUE with variance $t$. 
The extended Bessel kernel ${\bf K}_{J_\nu}$ 
is derived as the hard edge scaling limit \cite{FNH99,TW04}:
\begin{equation}
\mbK_N^{(\nu)} \left( N+s, x; N+t, y \right)
\to {\bf K}_{J_\nu}(s,x;t,y).
\quad N\to\infty,
\label{limit:hardedge}
\end{equation}
which implies
\begin{eqnarray}
&& (\{\Xi^{(\nu)}(t)\}_{ t \in [0,\infty)}, \P_{ \mu^{(\nu)}_{N,N} } )
\, \fd \, (\{\Xi^{(\nu)}(t)\}_{t\in [0,\infty)}, {\bf P}^{(\nu)}), 
\quad N \to \infty,
\nonumber\\
\label{conv:vaguely2}
&& \hskip 3cm 
\quad \mbox{in the vague topology}.
\end{eqnarray}

The extended Airy kernel
is derived as the soft-edge scaling limit \cite{PS02,Joh03,NKT03}:
\begin{eqnarray}
&&\mbK_N \Big( N^{1/3} +s, 2N^{2/3}+N^{1/3}s -\frac{s^2}{4}+x
; N^{1/3} +t, 2N^{2/3}+N^{1/3}t -\frac{t^2}{4}+y \Big)
\nonumber\\
&&\qquad\qquad\qquad\qquad\qquad\qquad
\to {\bf K}_{\Ai}(s,x;t,y),
\quad N\to\infty.
\label{lim:softedge}
\end{eqnarray}
We introduce the nonnegative function 
$f_{\rm sc}(x)= (2/\pi) \sqrt{1-x^2}{\bf 1}(|x|\le 1)$,
and put 
$$
\displaystyle{
\rho^N_{\rm sc}(t,x)
=\sqrt{\frac{N}{4t}}f_{\rm sc}\left(\frac{x}{2\sqrt{Nt}}\right) 
= \frac{1}{2\pi t}\sqrt{4tN-x^2}{\bf 1}(|x|\le 2\sqrt{tN})}, 
$$
associated with Wigner's semicircle law of eigenvalue 
distribution in the GUE \cite{Meh04}.
Then we put
\begin{eqnarray}
\widehat{\rho}^N_{\rm sc}(x) &=& \rho^N_{\rm sc}(N^{1/3},2N^{2/3}+x)
=\frac{N^{1/3}}{2}f_{\rm sc}\left(1+\frac{x}{2N^{2/3}}\right){\bf 1}(x<0)
\nonumber\\
&=& \frac{1}{\pi}\sqrt{-x \left(1+ \frac{x}{4N^{2/3}}\right)}{\bf 1}(x<0).
\label{def:circleN}
\end{eqnarray}
Here we note that
$\widehat{\rho}^N_{\rm sc}(x) \nearrow \widehat{\rho}(x)$, $N\to\infty$
and
\begin{equation}
\int_{\R} dx \ \widehat{\rho}^N_{\rm sc}(x) =N,
\quad \mbox{and} \quad
\int_{-4N^{2/3}}^0 dx \ \frac{\widehat{\rho}^N_{\rm sc}(x)}{-x} =N^{1/3}.
\nonumber
\end{equation}
Then (\ref{lim:softedge}) implies
\begin{eqnarray}
&& (\{\Xi_{\widehat{\rho}^N_{\rm sc}}(t)\}_{t \in [0,\infty)}, \P_{\tau_{-2N^{2/3}} \mu^{\rm GUE}_{N,N^{1/3}}})
\, \fd \,
(\{\Xi(t)\}_{t\in [0,\infty)}, {\bf P}_{{\rm Ai}}), \quad N \to \infty,
\nonumber\\
\label{conv:vaguely3}
&& \hskip 3cm  
\quad \mbox{in the vague topology}.
\end{eqnarray}

\SSC{Proofs of Results}\label{sec4}

\subsection{Proof of (i) of Theorems \ref{Theorem:Dyson} and 
\ref{Theorem:Bessel}}

To obtain (i) of the theorems 
it is enough to show $\mu_{\sin}(\mY)=1$ and $\mu_{J_\nu}(\mY^+)=1$,
since the other claims are derived from the facts 
that the operator $T_t$ and $T_t^{(\nu)}$ are of contraction and 
the set of all bounded continuous functions on $\mM$ 
is dense in $L^2(\mM,\mu_{\sin})$ and  $L^2(\mM,\mu_{J_\nu})$.

First we consider the probability measure $\mu_{\sin}$.
Since its density $\rho_{\sin}$ is constant, 
the condition (\ref{eqn:71_i}) is satisfied for $\rho=\rho_{\sin}$.
Note that any correlation $\rho_m$ is bounded, 
because the kernel ${K_{\sin}}$ is bounded.
Then if we take $\kappa \in (1/2,1)$ and $m\in\N$ satisfying $(1-\kappa)m>1$,
we see that the condition (\ref{eqn:73}) is satisfied 
for $\rho=\rho_{\sin}$.
From Lemma \ref{lemma:correlation_1} with $k=2$, 
we see that the conditions (\ref{eqn:74}) and (\ref{eqn:75}) 
are satisfied with $m'=4, p=2$ and $\rho=\rho_{\sin}$. 
Thus, $\mu_{\sin}(\mY)=1$.
For $\mu_{J_\nu}$, its density $\rho^{(\nu)}$ is 
bounded and satisfies 
$\int_0^\infty \{\rho^{(\nu)}(x)/x \} dx <\infty$. 
Hence, by the same argument as above we obtain $\mu_{J_\nu}(\mY^+)=1$.


\vskip 3mm
\noindent{\bf Remark 3.} \quad
When $\mu_\lambda$ is the Poisson point process 
with an intensity measure $\lambda dx$, $\lambda>0$,
$\rho_m(\x_m)=\lambda^m$. Then we can readily confirm that
all assumptions in Lemma \ref{lemma:correltion_4} hold
with $m \in \N, \kappa \in (1/2, 1)$
satisfying $(1-\kappa) m > 1$ and with
$m'=4$ and $p=2$.
Then $\mu_\lambda (\mY)=1$.
We can also show that measures such as 
Gibbs states with regular conditions are applicable
to Lemma \ref{lemma:correltion_4}.
\subsection{Proof of (ii) and (iii) of Theorems 
\ref{Theorem:Dyson} and \ref{Theorem:Bessel}}

The following lemma states stronger properties than 
(\ref{conv:vaguely1}) and (\ref{conv:vaguely2})
and is the key to proving (ii) and (iii) of
Theorems \ref{Theorem:Dyson} and \ref{Theorem:Bessel}.
The proof of this lemma is given in the next subsection.

\begin{lem}
\label{lem:KEY_LEMMA}
{\rm (i)}
For any $0 \equiv t_0 < t_1 < t_2 < \cdots < t_M < \infty$, $M \in \N_0$,
and bounded functions $g_j, 0 \leq j \leq M$, which is 
$\Phi_0$-moderately continuous on $\mY_{m,\kappa}$
for any $m\in\N$ and $\kappa\in (1/2,1)$,
\begin{equation}
\lim_{N\to\infty} \E_{\mu^{\rm GUE}_{N,2N/\pi^2}} 
\left[ \prod_{j=0}^{M} g_j(\Xi(t_j)) \right]
= {\bf E}_{\sin} \left[ \prod_{j=0}^{M} g_j(\Xi(t_j)) \right].
\label{eqn:Phi-cont1}
\end{equation}
In particular, $(\{\Xi(t)\}_{t \in [0,\infty)}, \P_{\mu^{\rm GUE}_{N,2N/\pi^2}})
\fd (\{\Xi(t)\}_{t\in [0,\infty)}, {\bf P}_{\sin})$,
$N \to \infty$, in the sense of finite dimensional distributions in the $\Phi_0$-moderate topology.

\vskip 3mm 

\noindent {\rm (ii)} 
For any $0 \equiv t_0 < t_1 < t_2 < \cdots < t_M < \infty$, $M \in \N_0$,
and bounded functions $g_j, 0 \leq j \leq M$, which is 
$\Phi_0$-moderately continuous on $\mY^{+}_{m,\kappa}$
for any $m\in\N$ and $\kappa\in (1,2)$,
\begin{equation}
\lim_{N\to\infty} \E_{\mu^{(\nu)}_{N,N}} 
\left[ \prod_{j=0}^{M} g_j(\Xi(t_j)) \right]
= {\bf E}_{J_\nu} \left[ \prod_{j=0}^{M} g_j(\Xi(t_j)) \right].
\label{eqn:Phi-cont2}
\end{equation}
In particular, 
$( \{\Xi^{(\nu)}(t)\}_{t \in [0,\infty)}, \P_{\mu^{(\nu)}_{N,N}}) \fd 
(\{\Xi(t)\}_{t\in [0,\infty)}, {\bf P}_{J_\nu})$, 
$N \to \infty$, in the sense of finite dimensional distributions in the $\Phi_0$-moderate topology.
\end{lem}
\vskip 3mm

\noindent {\it Proof of (ii) and (iii) in Theorems \ref{Theorem:Dyson} 
and \ref{Theorem:Bessel}.} \ 
We first give the proof for Theorem \ref{Theorem:Dyson}.
Let $f_j, 0\le j \le M$ be bounded vaguely continuous functions on $\mM$
and $0=t_0 \le t_1<t_2<\cdots<t_M <\infty$, $M\in\N$.
Theorem \ref{Theorem:Dyson} (ii) is 
concluded from the equality
\begin{equation}
\int_{\mM} \mu_{\sin}(d\xi)
\E_{\xi} \left[\prod_{j=0}^M f_j(\Xi(t_j)) \right]
= \lim_{N \to \infty} 
\int_{\mM} \mu^{\rm GUE}_{N,2N/\pi^2}(d\xi)
\E_{\xi} \left[\prod_{j=0}^M f_j(\Xi(t_j)) \right].
\label{eqn:key_ii}
\end{equation}
Since $\E_{\xi} \left[\displaystyle{\prod_{j=0}^M} f_j(\Xi(t_j)) \right]$ 
is $\Phi_0$-moderately continuous on $\mY_{m,\kappa}$ 
for any $m\in\N$ and $\kappa\in (1/2,1)$ 
from Proposition \ref{Theorem:Dyson_CMP} (ii),
(\ref{eqn:key_ii}) is guaranteed by 
(\ref{eqn:Phi-cont1}) of Lemma \ref{lem:KEY_LEMMA} with $M=0$
and $g_0(\xi)= \E_{\xi} \left[\displaystyle{\prod_{j=0}^M} 
f_j(\Xi(t_j)) \right]$.
Then Theorem \ref{Theorem:Dyson} (ii) is proved.
Now we show Theorem \ref{Theorem:Dyson} (iii), 
which is derived from the relations:
\begin{eqnarray}
&&{\bf E}_{\sin} 
\Big[f_0(\Xi(0))f_1(\Xi(t_1))f_2(\Xi(t_2))
\cdots f_M(\Xi(t_M) \Big]
\nonumber\\
&&\quad = {\bf E}_{\sin} \Big[f_0(\Xi(0))f_1(\Xi(t_1))
\E_{\Xi(t_1)}\Big[f_2(\Xi(t_2-t_1)) 
\nonumber\\
&&\qquad\qquad\qquad
\cdots
\E_{\Xi(t_{M-1}-t_{M-2})}\Big[f_M(\Xi(t_M-t_{M-1}))
\Big]\cdots \Big]\Big]
\nonumber\\
&&\quad
=\langle f_0, T_{t_1}( f_1T_{t_2-t_1}(f_2\cdots T_{t_{M}-t_{M-1}}f_M)\cdots) 
\rangle_{\mu_{\sin}}.
\label{eqn:Markovian}
\end{eqnarray}
We show the first equality of (\ref{eqn:Markovian}) 
by induction with respect to $M$.
First we consider the case $M=2$.
By the Markov property of 
$(\{\Xi(t)\}_{t\in [0,\infty)}, \P_{\mu^{\rm GUE}_{N,2N/\pi^2}})$,
we have
\begin{eqnarray}
&&\E_{\mu^{\rm GUE}_{N,2N/\pi^2}} 
\Big[f_0(\Xi(0))f_1(\Xi(t_1))f_2(\Xi(t_2)) \Big]
\nonumber\\
&&\qquad= \E_{\mu^{\rm GUE}_{N,2N/\pi^2}} 
\Big[f_0(\Xi(0))f_1(\Xi(t_1))
\E_{\Xi(t_1)}\Big[f_2(\Xi(t_2-t_1)) \Big]\Big]
\nonumber
\end{eqnarray}
for $0\le t_1<t_2<\infty$.
Since
$$
\lim_{N\to\infty}\E_{\mu^{\rm GUE}_{N,2N/\pi^2}} 
\Big[f_0(\Xi(0))f_1(\Xi(t_1))f_2(\Xi(t_2)) \Big]
= {\bf E}_{\sin} \Big[f_0(\Xi(0))f_1(\Xi(t_1))f_2(\Xi(t_2)) \Big]
$$
by (\ref{conv:vaguely1}), it is enough to show
\begin{eqnarray}
&&\lim_{N\to\infty}\E_{\mu^{\rm GUE}_{N,2N/\pi^2}} \Big[f_0(\Xi(0))f_1(\Xi(t_1))\E_{\Xi(t_1)}\Big[f_2(\Xi(t_2-t_1)) \Big]\Big]
\nonumber\\
&& \qquad 
= {\bf E}_{\sin} \Big[f_0(\Xi(0))f_1(\Xi(t_1))\E_{\Xi(t_1)}\Big[f_2(\Xi(t_2-t_1)) \Big]\Big]
\label{convergence:Markov}
\end{eqnarray}
for the proof. Since 
$\E_{\xi}\Big[f_2(\Xi(t_2-t_1)) \Big]=T_{t_2-t_1}f_2(\xi)$
is $\Phi_0$-moderately continuous on $\mY_{m,\kappa}$
for any $m\in\N$ and $\kappa\in (1/2,1)$, (\ref{eqn:Phi-cont1}) 
of Lemma \ref{lem:KEY_LEMMA} with $M=1$ 
and $g_0=f_0, g_1=f_1T_{t_2-t_1}f_2$ gives (\ref{convergence:Markov}). 
Thus, we have obtained the case $M=2$.
Next we suppose that the first equality of (\ref{eqn:Markovian}) 
is satisfied in case $M=k\in\N$.
By using the same argument as in the case $M=2$,
in which (\ref{eqn:Phi-cont1}) of 
Lemma \ref{lem:KEY_LEMMA} is used in this case with $M=1$ and
$g_0=\prod_{j=0}^{k-1} f_j, 
g_1=f_k T_{t_{k+1}-t_k} f_{k+1}$, we have
$$
{\bf E}_{\sin} \left[\prod_{j=0}^{k+1} f_j(\Xi(t_j)) \right]
={\bf E}_{\sin} \left[\prod_{j=0}^kf_j(\Xi(t_j))
\E_{\Xi(t_k)}\Big[f_{k+1}(\Xi(t_{k+1}-t_k)) \Big] \right]
$$
with $t_0=0$.
From the assumption of the induction the right hand side of
the above equality equals to 
\begin{eqnarray}
&&{\bf E}_{\sin} \Big[f_0(\Xi(0))f_1(\Xi(t_1))
\E_{\Xi(t_1)}\Big[f_2(\Xi(t_2-t_1)) 
\nonumber\\
&&\qquad\qquad\qquad
\cdots
\E_{\Xi(t_{k}-t_{k-1})}\Big[f_{k+1}(\Xi(t_{k+1}-t_{k}))
\Big]\cdots \Big]\Big].
\nonumber
\end{eqnarray}
Then we obtain the first equality of 
(\ref{eqn:Markovian}) in case $M=k+1$, and
the induction is completed.

The second equality of (\ref{eqn:Markovian}) is derived
from the definition (\ref{eqn:defTt}) of the
operator $T_t$ and (\ref{eqn:inner_product}).
That is, 
\begin{eqnarray}
&& {\bf E}_{\sin} \Big[f_0(\Xi(0))f_1(\Xi(t_1))
\E_{\Xi(t_1)}\Big[f_2(\Xi(t_2-t_1)) 
\nonumber\\
&&\qquad\qquad\qquad
\cdots
\E_{\Xi(t_{M-1}-t_{M-2})}\Big[f_M(\Xi(t_M-t_{M-1}))
\Big]\cdots \Big]\Big]
\nonumber\\
&&\quad
= {\bf E}_{\sin} \Big[
f_0(\Xi(0)) f_1(\Xi(t_1)) 
T_{t_2-t_1}(f_2\cdots T_{t_{M}-t_{M-1}}f_M)\cdots)(\Xi(t_1)) \Big]
\nonumber\\
&&\quad
=\langle f_0, T_{t_1}( f_1T_{t_2-t_1}
(f_2\cdots T_{t_{M}-t_{M-1}}f_M)\cdots) 
\rangle_{\mu_{\sin}}.
\nonumber
\end{eqnarray}
The proof of Theorem \ref{Theorem:Dyson}
is completed.

The proof of Theorem \ref{Theorem:Bessel} is obtained by the same argument,
in which (\ref{conv:vaguely2}) and (\ref{eqn:Phi-cont2}) should be used
instead of (\ref{conv:vaguely1}) and  (\ref{eqn:Phi-cont1}).
\qed

\subsection{Proof of Lemma \ref{lem:KEY_LEMMA}}
\vskip 3mm

We introduce subsets of $\mY$,
$$
\mY^{\gamma,L_0}_{\kappa, m}=\bigg\{ \xi\in\mM : m(\xi,\kappa) \le m, \
\bigg|\int_{|x|\ge L}\frac{\xi(dx)}{x}\bigg|\le L^{-\gamma}, 
\ L\ge L_0 \bigg\},
$$
with 
$\kappa\in(1/2,1)$ $\gamma>0$, $m,L_0\in\N$.
Then we can prove the following lemma.

\begin{lem} 
\label{lemma:4_5_4}

{\rm (i)} For any $t>0$, 
$\displaystyle{\lim_{m\to\infty}\lim_{L_0\to\infty} 
\min_{N\in\N}\mu^{\rm GUE}_{N,2N/\pi^2+t} 
\Big( \mY^{\gamma,L_0}_{\kappa, m} \Big)=1}$
for some $\kappa \in (1/2,1)$ and $\gamma>0$.

{\rm (ii)} For any $t>0$, 
$\displaystyle{\lim_{m\to\infty} \min_{N\in\N}\mu^{(\nu)}_{N,N+t} 
\Big( \mY^{+}_{\kappa, m} \Big)=1}$
for some $\kappa \in (1,2)$.

\end{lem}

\noindent {\it Proof.}
\quad 
Here we give the proof of (i) only, since
(ii) will be proved by the similar argument to the latter half
of the proof of (i).
For any fixed $t>0$ we put
\begin{equation}
\rho^N_{\rm GUE}(t,x) \equiv 
\mbK_N(t,x;t,x)= \frac{1}{\sqrt{2t}}\sum_{k=0}^{N-1}
\varphi_{k}\left(\frac{x}{\sqrt{2t}}\right)^2,
\label{def:rhoGUE}
\end{equation}
and
$$
\rho^N(t,x) \equiv \rho_{\rm GUE}^N \left(\frac{2N}{\pi^2}+t, x \right),
$$
Then $\rho^N(t, x)$ is a symmetric function of 
$x$ and bounded with respect to $N$ and $x$.
Since $\mu^{\rm GUE}_{N,2N/\pi^2+t}$ is a determinantal point process, 
by Lemma \ref{lemma:correlation_1} we have
$$
\int_{\mM}\mu^{\rm GUE}_{N,2N/\pi^2+t}(d\xi)
\left|\xi([0,L))-\int_0^L \rho^N(t;x)dx \right|^4
\le CL^2
$$
with a positive constant $C$, which is independent of $N$.
By Chebyshev's inequality we have
\begin{equation}
\label{estimate:muN}
\mu^{\rm GUE}_{N,2N/\pi^2+t} \left(
\left|\xi([0,L))-\int_0^L \rho^N(t;x)dx \right|\ge L^{7/8}
\right)
\le CL^{-3/2},
\end{equation}
and so
$$
\mu^{\rm GUE}_{N,2N/\pi^2+t}
\left( \left|\xi([0,L))-\int_0^L \rho^N(t;x)dx \right|
\le L^{7/8},
\ \forall L\ge L_0 \right)
\ge 1- C' L_0^{-1/2}.
$$
Similarly, we have
$$
\mu^{\rm GUE}_{N,2N/\pi^2+t}
\left( \left|\xi((-L,0])-\int_{-L}^0 \rho^N(t;x)dx \right|
\le L^{7/8},
\ \forall L\ge L_0 \right)
\ge 1- C' L_0^{-1/2}.
$$
Note that 
$$
\int_{L\le |x|\le L'}  \frac{\rho^N(t;x)dx}{x}=0, \quad L'>L,
$$
from the fact that $\rho^N(t, x)$ is symmetric in $x$. 
Then, by the same procedure in the proof 
of Lemma \ref{lemma:correlation_3} with $C_1=1$ 
and $\varepsilon=7/8$, we have
\begin{equation}
\label{estimate:muNt_1}
\mu^{\rm GUE}_{N,2N/\pi^2+t}
\left( \left|\int_{|x|\ge L}  \frac{\xi(dx)}{x} \right|
\le 24L^{-1/8}, \ \forall L\ge L_0 \right)
\ge 1- 2C' L_0^{-1/2}.
\end{equation}

By Chebyshev's inequality with Lemma \ref{lemma:correlation_1},
for any $p\in\N$
\begin{eqnarray}
\mu^{\rm GUE}_{N,2N/\pi^2+t}
\bigg(\xi \Big([g^{\kappa}(k),g^{\kappa}(k+1)] \Big)\ge m\bigg)
&\le& 
\frac{\big(3 \rho^N(t, [g^{\kappa}(k),g^{\kappa}(k+1)])\big)^p}
{|m- \rho^N(t, [g^{\kappa}(k),g^{\kappa}(k+1)])|^{2p}}
\nonumber\\
&\le&Ck^{(\kappa -1)p}m^{-2p},
\label{estimate:gkapp}
\end{eqnarray}
where $\rho^N(t,D) \equiv \int_{D}dx \rho^N(t, x)$ 
and $C$ is a constant independent of $N$ and $k$.
If $p$ is large enough to satisfy $(1-\kappa)p >1$, we have
$$
\mu^{\rm GUE}_{N,2N/\pi^2+t}\bigg(\max_{k\in\Z}
\xi \Big([g^{\kappa}(k),g^{\kappa}(k+1)] \Big) \ge m \bigg)\le C'm^{-2p},
$$
and so we see
\begin{equation}
\label{estimate:muNt_2}
\lim_{m\to\infty}\min_{N\in\N}\mu^{\rm GUE}_{N,2N/\pi^2+t}\bigg(
\max_{k\in\Z}
\xi \Big([g^{\kappa}(k),g^{\kappa}(k+1)] \Big) \le m
\bigg)= 1.
\end{equation}
Combining the above estimates (\ref{estimate:muNt_1}) 
and (\ref{estimate:muNt_2}),
we obtain (i) of the lemma.
\qed
\vskip 3mm

\noindent {\it Proof of Lemma \ref{lem:KEY_LEMMA}.}\quad 
Since the proofs of (i) and (ii) are similar, here 
we only give the proof of (i). 
Let $d(\cdot,\cdot)$ be a metric on $\mM$ associated 
with the vague topology. First remind that $\xi_n$ converges 
$\Phi_0$-moderately to $\xi$ if the conditions (\ref{eqn:THC3}) 
and (\ref{eqn:THC4}) are satisfied. Then from the definition of 
$\mY^{\gamma,L_0}_{\kappa, m}$, we see that for any $\varepsilon >0$ 
there exists $\delta >0$ such that 
$$
|g_j(\xi)-g_j(\eta)|<\varepsilon, \quad 0 \leq j \leq M,
$$
for any $\xi, \eta \in \mY^{\gamma,L_0}_{\kappa, m}$ with $d(\xi, \eta)<\delta$.
Here we use the fact that the closure of $\{\xi \in \mM : m(\xi, \kappa) \leq m\}$ is a compact subset of $\mM$ to ensure that $\delta$ does not depend on $\xi$ or $\eta$. Then, from the fact (\ref{conv:vaguely1}), we can show that
\begin{eqnarray}
&& \left| 
\lim_{N\to\infty} \E_{\mu^{\rm GUE}_{N,2N/\pi^2}} 
\left[ \prod_{j=0}^{M} g_j(\Xi(t_j)) \right] - {\bf E}_{\sin} \left[ \prod_{j=0}^{M} g_j(\Xi(t_j)) \right]  \right|
\nonumber\\
&& \le \prod_{j=0}^M \sup_{\xi}|g_j(\xi)|
\left\{ (M+1) \mu_{\sin}(\mM \setminus \mY^{\gamma,L_0}_{\kappa, m}) + \sum_{j=0}^{M} \max_{N\in\N} \mu^{\rm GUE}_{N,2N/\pi^2+t_j} (\mM \setminus \mY^{\gamma,L_0}_{\kappa, m}) \right\}
\nonumber
\end{eqnarray}
by using the Skorohod representation theorem,
which can be applied to distributions on Polish spaces \cite{BD83}.
Hence, by Lemma \ref{lemma:4_5_4} we obtain the lemma.
\qed

\subsection{Proof of Theorem \ref{Theorem:Airy}}

We put
\begin{equation}
\quad
\rho_{\cal A}^N (x) =\rho_{\rm GUE}^N(N^{1/3},2N^{2/3}+x),
\label{def:rho_ca}
\end{equation}
where $\rho_{\rm GUE}^N(t,x)$ is defined in (\ref{def:rhoGUE}).
The soft-edge scaling limit (\ref{lim:softedge}) implies that
$$
\displaystyle{\lim_{N\to\infty}\rho_{\cal A}^N (x)
=\rho_{\rm Ai}(x) \equiv K_{\rm Ai}(x,x)}.
$$
We also see the following estimates, whose proof will be
given in Section \ref{proof:L32}.

\begin{lem}
\label{lemma:correlation_2}
There exists a positive constant $C_3$ such that
\begin{equation}
|\rho_{\rm Ai}(x)- \widehat{\rho}(x)|\le  \frac{C_3}{|x|},
\quad x\in \R,
\label{asym:Airy}
\end{equation}
and
\begin{equation}
|\rho_{\cal A}^N(x) -\widehat{\rho}_{\rm sc}^N(x)|\le \frac{C_3}{|x|},
\quad x\in \R, \ N\in\N.
\label{asym:AN}
\end{equation}
In particular
\begin{equation}
\max_{N\in\N} \int_{|x| \ge L} dx \,
\left| \frac{\rho^N_{\cal A}(x) -\widehat{\rho}_{\rm sc}^N(x)}{x}\right|
\le \frac{2C_3}{L}.
\label{eqn:dens2}
\end{equation}
\end{lem}

For the proof of (i) we first show  $\mu_{{\rm Ai}}(\mY^{\cal A})=1$.
Let $\rho_{\Ai}$ be the density function of $\mu_{\Ai}$.
Note that $\rho_{\Ai}(x)={\cal O}(|x|^{1/2})$, 
$x\to -\infty$ from (\ref{asym:Airy}). 
Using Lemma \ref{lemma:correlation_1} with $k=3$,
we see that (\ref{eqn:74}) and (\ref{eqn:75})
are satisfied with $m'=6, p=4.5$.
Note the correlation inequality $\rho_m(\x_m)\le \prod_{j=1}^m \rho(x_j)$
(see Proposition 4.3 in \cite{ST03}).
If we take $\kappa \in (1/2,2/3)$ 
and $m\in\N$ satisfying $(3\kappa/2-1)m<-1$,
we see that the condition (\ref{eqn:73}) is satisfied for
the correlation functions $\rho_m$ of $\mu_{\Ai}$.
The condition (\ref{eqn:71_ii}) 
with $\overline{\rho}=\widehat{\rho}$ 
is derived from (\ref{asym:Airy}) immediately.
Hence, from Lemma \ref{lemma:correltion_4} (iii) 
we obtain the desired result.

\vskip 3mm

For the proof of (ii) and (iii) we introduce subsets of 
$\mY^{\cal A}$,
$$
\mY^{{\cal A},\gamma,L_0}_{\kappa,m}=
\left\{ \xi\in\mM : m(\xi,\kappa) \le m, \
\left|\int_{|x|\ge L}\frac{\widehat{\rho}_{\rm sc}^{\xi(\R)}(x)dx 
-\xi(dx)}{x}\right|\le L^{-\gamma},  \ L\ge L_0 \right\}
$$
with $\kappa\in(1/2,2/3]$ $\gamma>0$, $m,L_0\in\N$.
Then the following lemma is established.

\begin{lem} 
\label{lemma:4_5_4_A} We have
$$
\lim_{m\to\infty}\lim_{L_0\to\infty} 
\min_{N\in\N}\mu^{\rm GUE}_{N,N^{1/3}} 
\Big( \tau_{2N^{2/3}}\mY^{{\cal A},\gamma,L_0}_{\kappa,m} \Big)=1
$$
for some $\kappa \in (1/2,2/3)$ and $\gamma>0$.
\end{lem}

\noindent {\it Proof.}
Note that $\int_0^\infty \rho_{\cal A}^N(x) dx <\infty$ and
$\int_{-L}^0 \rho_{\cal A}^N(x) dx \le CL^{3/2}$ 
from Lemma \ref{lemma:correlation_2}.
Since $\mu^{\rm GUE}_{N,N^{1/3}}$ is 
a determinantal point process, by Lemma \ref{lemma:correlation_1}
$$
\int_{\mM}\tau_{-2N^{2/3}}\mu^{\rm GUE}_{N,N^{1/3}}(d\xi)
\bigg|\xi(-L,\infty))-\int_{-L}^\infty \rho_{\cal A}^N(x)dx\bigg|^6
\le CL^{9/2}
$$
with a positive constant $C$, which is independent of $N$.
By Chebyshev's inequality we have
\begin{equation}
\label{estimate:muNA}
\tau_{-2N^{2/3}}\mu^{\rm GUE}_{N,N^{1/3}} \left(
\left|\xi((-L,\infty))-\int_{-L}^\infty\rho_{\cal A}^N(x)dx \right|
\ge L^{23/24}\right)
\le \frac{C}{L^{5/4}},
\end{equation}
and so
$$
\tau_{-2N^{2/3}}\mu^{\rm GUE}_{N,N^{1/3}}
\left( \left|\xi((-L,\infty))-\int_{-L}^\infty \rho_{\cal A}^N(x)dx \right|
\le L^{23/24},
\ \forall L\ge L_0 \right)
\ge 1- \frac{C'}{ L_0^{1/4}}.
$$
By using the estimate (\ref{asym:AN}) and Lemma 
\ref{lemma:correlation_3} (ii) with $C_1=1$ 
and $\varepsilon=23/24$, we see that
\begin{equation}
\label{estimate:muNt_1_A}
\tau_{-2N^{2/3}}\mu^{\rm GUE}_{N,N^{1/3}}
\left( 
\left|\int_{|x|\ge L}  
\frac{\widehat{\rho}^N_{\rm sc}(x)dx -\xi(dx)}{x} \right|
\le \frac{C"}{L^{1/24}}, \ \forall L\ge L_0 \right)
\ge 1- \frac{2C'}{ L_0^{1/4}}
\end{equation}
with $C''=72+C_3$.
Noting the estimate (\ref{asym:AN}), 
by the same argument deriving (\ref{estimate:gkapp}) we will obtain
$$
\tau_{-2N^{2/3}}\mu^{\rm GUE}_{N,N^{1/3}}
\bigg(\xi \Big([g^{\kappa}(k),g^{\kappa}(k+1)] \Big)\ge m\bigg)
 \le Ck^{(3\kappa/2 -1)p}m^{-2p}.
$$
and
\begin{equation}
\label{estimate:muNt_2_A}
\lim_{m\to\infty}\min_{N\in\N}
\tau_{-2N^{2/3}}\mu^{\rm GUE}_{N,N^{1/3}}
\bigg( \max_{k\in\Z} \xi \Big([g^{\kappa}(k),g^{\kappa}(k+1)] \Big) \le m
\bigg)= 1.
\end{equation}
Combining the above estimates (\ref{estimate:muNt_1_A}) 
and (\ref{estimate:muNt_2_A}),
we obtain the lemma.
\qed
\vskip 3mm

\noindent {\it Proof of (ii) and (iii) of Theorem \ref{Theorem:Airy}.}\quad 
For any fixed $t>0$ we put
\begin{eqnarray}
&&\widehat{\rho}_{\rm sc}^N(t,x)
= \sqrt{\frac{N^{1/3}}{N^{1/3}+t}}\widehat{\rho}_{\rm sc}^N
\left(
\sqrt{\frac{N^{1/3}}{N^{1/3}+t}} \ x
\right),
\nonumber\\
&&\rho_{\cal A}^N(t,x)=\rho_{\rm GUE}^N \left(N^{1/3}
+t, 2N^{1/2}(N^{1/3}+t)^{1/2}+x \right), 
\quad N \in \N,
\nonumber
\end{eqnarray}
and
$$
\mY^{{\cal A},\gamma,L_0}_{\kappa,m}(t)
=\left\{ \xi\in\mM : m(\xi,\kappa) \le m, \
\left|\int_{|x|\ge L}\frac{\widehat{\rho}_{\rm sc}^{\xi(\R)}(t,x)dx 
-\xi(dx)}{x}\right|\le L^{-\gamma},  \ L\ge L_0 \right\}.
$$
Note that $\widehat{\rho}^N_{\rm sc}(t, \cdot)$ 
is a nonnegative function such that 
$\int_{\R} dx \, \widehat{\rho}^N_{\rm sc}(t, x)=N$,
$\widehat{\rho}^N_{\rm sc}(t, x) \nearrow \widehat{\rho}(x)$, $N\to\infty$.

By using the scaling property of the Dyson model, 
from Lemma \ref{lemma:4_5_4_A}
we have
$$
\lim_{m\to\infty}\lim_{L_0\to\infty} 
\min_{N\in\N}\mu^{\rm GUE}_{N,N^{1/3}+t} 
\Big( \tau_{2N^{1/2}(N^{1/3}+t)^{1/2}}
\mY^{{\cal A}, \gamma, L_0}_{\kappa, m}(t) \Big)=1,
$$
for some $\kappa \in (1/2,2/3)$ and $\gamma>0$.
Note that $2N^{1/2}(N^{1/3}+t)^{1/2}-2N^{2/3}+N^{1/3}t -t^2/4 = {\cal O}(N^{-1/3})$. Since 
$$
\left|\int_{|x|\ge L}\frac{\widehat{\rho}_{\rm sc}^{N}(t,x)dx 
-\widehat{\rho}_{\rm sc}^{N}(t,x+\varepsilon)dx}{x}\right|
\le \varepsilon L^{-1/2},
$$
we have
$$
\lim_{m\to\infty}\lim_{L_0\to\infty} 
\min_{N\in\N}\mu^{\rm GUE}_{N,N^{1/3}+t} 
\Big( \tau_{ 2N^{2/3}-tN^{1/3}+t^2/4}
\mY^{{\cal A}, \gamma, L_0}_{\kappa, m}(t) \Big)=1,
$$
and so
\begin{equation}
\lim_{m\to\infty}\lim_{L_0\to\infty} \min_{N\in\N}
\P_{\tau_{-2N^{2/3}}\mu^{\rm GUE}_{N,N^{1/3}}}
\left(
\Xi_{\widehat{\rho}_{\rm sc}^N}(t) \in\mY^{{\cal A},\gamma,L_0}_{\kappa,m}(t)
\right)=1
\label{llm_AIry}
\end{equation}
for some $\kappa \in (1/2,2/3)$ and $\gamma>0$.

Let $0 \equiv t_0 < t_1 < t_2 < \cdots < t_M < \infty$, $M \in \N_0$.
Suppose that $g_j, 0 \leq j \leq M$ are bounded functions on $\mM$ such that
$g_j(\xi^N) \to g_j(\xi)$, $N\to\infty$, 
if $\xi\in \mY^{\cal A}$ and $\xi^N \in\mM$ with $\xi^N(\R)=N$
satisfy that 
$\xi^N\to \xi$ vaguely,  
(\ref{condition:convergence}) holds
with $\widehat{\rho}^N(\cdot)=\widehat{\rho}^N_{\rm sc}(t_j,\cdot)$, 
and $\max_{N\in\N} m(\xi^N,\kappa)\le m$ for some $m\in\N$ and $\kappa\in (1/2,1)$.
From (\ref{llm_AIry}),
we can show 
\begin{equation}
\lim_{N\to\infty} \E_{\tau_{-2N^{2/3}}\mu^{\rm GUE}_{N,N^{1/3}}} 
\left[ \prod_{j=0}^{M} g_j(\Xi(t_j)) \right]
= {\bf E}_{\rm A_i} \left[ \prod_{j=0}^{M} g_j(\Xi(t_j)) \right].
\nonumber
\end{equation}
by the same argument as the proof of Lemma \ref{lem:KEY_LEMMA}.
Therefore, by applying (ii) of Proposition \ref{Theorem:JSP}
we can show (ii) and (iii) of Theorem \ref{Theorem:Airy}
by the same procedure as in the proof of Theorem \ref{Theorem:Dyson}.
\qed

\SSC{Proof of Lemma \ref{lemma:correlation_2}} \label{proof:L32}

\subsection{Proof of (\ref{asym:Airy})}

We use the following asymptotic expansions of the Airy functions 
in the classical sense of Poincar\'e \cite{VS04,Olv54a,Olv54b}: 
For $x\gg 1$
\begin{eqnarray}
&&{\rm Ai}(x)\approx \frac{e^{-\frac{2}{3}x^{3/2}}}{2\pi^{1/2}x^{1/4}}
L\left(-\frac{2}{3}x^{3/2}\right),
\quad
{\rm Ai}'(x)\approx -\frac{x^{1/4}e^{-\frac{2}{3}x^{3/2}}}{2\pi^{1/2}}
M\left(-\frac{2}{3}x^{3/2}\right),
\nonumber\\
&&{\rm Ai}(-x)\approx\frac{1}{\pi^{1/2}x^{1/4}}
\left[
\sin \left(\frac{2}{3}x^{3/2}-\frac{\pi}{4}\right)Q
\left(\frac{2}{3}x^{3/2}\right) 
+ \cos\left(\frac{2}{3}x^{3/2}-\frac{\pi}{4}\right)
P\left(\frac{2}{3}x^{3/2}\right) \right],
\nonumber\\
&&{\rm Ai}'(-x)\approx\frac{x^{1/4}}{\pi^{1/2}}
\left[
\sin \left(\frac{2}{3}x^{3/2}-\frac{\pi}{4}\right)R
\left(\frac{2}{3}x^{3/2}\right) 
- \cos\left(\frac{2}{3}x^{3/2}
-\frac{\pi}{4}\right)S\left(\frac{2}{3}x^{3/2}\right) \right],
\nonumber
\end{eqnarray}
where 
$L(z), M(z), P(z), Q(z), R(z)$ and $S(z)$ are functions defined by
\begin{eqnarray}
&&L(z)=\sum_{k=0}^\infty \frac{u_{k}}{z^{k}}
=1+\frac{5}{72 z} 
+ {\cal O}\Big(\frac{1}{z^2}\Big), \
M(z)=\sum_{k=0}^\infty \frac{v_{k}}{z^{k}}
=1- \frac{7}{72z}
+ {\cal O}\Big(\frac{1}{z^2}\Big),
\nonumber\\
&&P(z)=\sum_{k=0}^\infty (-1)^k \frac{u_{2k}}{z^{2k}}
=1+ {\cal O}\Big(\frac{1}{z^2}\Big), \,
Q(z)=\sum_{k=0}^\infty (-1)^k \frac{u_{2k+1}}{z^{2k+1}}
=\frac{5}{ 72z} 
+ {\cal O}\Big(\frac{1}{z^3}\Big),
\nonumber\\
&&R(z)=\sum_{k=0}^\infty (-1)^k \frac{v_{2k}}{z^{2k}}
=1+ {\cal O}\Big(\frac{1}{z^2}\Big), \,
S(z)=\sum_{k=0}^\infty (-1)^k \frac{v_{2k+1}}{z^{2k+1}}
=-\frac{7}{72z} 
+ {\cal O}\Big(\frac{1}{z^3}\Big),
\nonumber
\end{eqnarray}
with $\displaystyle{u_k=\frac{\Gamma(3k+1/2)}{54^k k! \Gamma(k+1/2)}}$,
$\displaystyle{v_k=-\frac{6k+1}{6k-1}u_k}$.
Then
\begin{eqnarray}
\rho_{\rm Ai}(x)&=& ({\rm Ai}'(x))^2 - x({\rm Ai}(x))^2
={\cal O}\left(e^{-\frac{4}{3}x^{3/2}} \right), \quad x\to\infty,
\nonumber
\end{eqnarray}
and
\begin{eqnarray}
&&\rho_{\rm Ai}(-x)= ({\rm Ai}'(-x))^2 + x ({\rm Ai}(-x))^2
\nonumber\\
&& =\frac{\sqrt{x}}{\pi}
\Bigg[
\bigg\{
\sin \left(\frac{2}{3}x^{3/2}-\frac{\pi}{4}\right)
\bigg( 1+ {\cal O}\bigg(\frac{1}{x^3}\bigg) \bigg)
- \cos\left(\frac{2}{3}x^{3/2}-\frac{\pi}{4}\right)
\bigg(-\frac{7}{72x^{3/2}} + {\cal O}\bigg(\frac{1}{x^{9/2}}\bigg)
\bigg)
\bigg\}^2
\nonumber\\
&&+\bigg\{
\sin \left(\frac{2}{3}x^{3/2}-\frac{\pi}{4}\right)
\bigg(\frac{5}{72x^{3/2}} 
+ {\cal O}\bigg(\frac{1}{x^{9/2}}\bigg)
\bigg)
+ \cos\left(\frac{2}{3}x^{3/2}-\frac{\pi}{4}\right)
\bigg( 1+ {\cal O}\bigg(\frac{1}{x^3}\bigg) \bigg)
\bigg\}^2
\Bigg]
\nonumber\\
&& =\frac{\sqrt{x}}{\pi}
\Bigg[
1+ \sin \left(\frac{2}{3}x^{3/2}-\frac{\pi}{4}\right)\cos \left(\frac{2}{3}x^{3/2}-\frac{\pi}{4}\right)
\frac{1}{3x^{3/2}}
+ {\cal O}\bigg(\frac{1}{x^3}\bigg)
\Bigg], \quad x \to \infty.
\nonumber
\end{eqnarray}
Hence, we obtain
\begin{equation}
\bigg|\rho_{\rm Ai}(x) - \widehat{\rho}(x)\bigg|
= {\cal O}\left(\frac{1}{|x|}\right),
\quad x\to\infty.
\nonumber
\end{equation}


\subsection{Proof of (\ref{asym:AN})}

For proving (\ref{asym:AN}) we use the asymptotic behavior of 
Hermite polynomials given in Plancherel and Rotach \cite{PR29}, 
in which the Hermite polynomials are defined as 
$\widehat{H}_n(x)=(-1)^n e^{x^2/2}
(d^n/dx^n) e^{-x^2/2}$, 
whereas our definition is
$H_n(x)=(-1)^n e^{x^2} (d^n/dx^n)e^{-x^2}$. 
We should note the relation that 
$\widehat{H}_n(x)=2^{-n/2}H_n\left(x/\sqrt{2} \right)$.
We introduce the following polynomials 
$\phi_{n}(z)$ and $\psi_{n p}(z)$
determined by the expansions 
\begin{eqnarray}
&&
\exp\Big[ 
z \sum_{m=3}^\infty \frac{(-1)^m}{m}\tau^{m-2}
\Big]
=\sum_{n=0}^\infty \phi_n (z)\tau^n, \quad
\nonumber
\\
&&
\left(\sum_{k =1}^\infty \frac{1}{k}\tau^{k-1}\right)^p
\exp \Big[
z \sum_{m=4}^\infty \frac{1}{m} \tau^{m-3}
\Big]
=\sum_{n=0}^\infty \psi_{n p}(z)\tau^n,
\nonumber
\end{eqnarray}
for  $|\tau|<1$. For example, 
$\phi_0(z)=1$, $\phi_1(z)=-z/3$, 
$\phi_2(z)=z/4+z^2/18$, 
$\phi_3(z)=-z/5-z^2/12-z^3/162$,
$\psi_{0 p}(z)=1$, $\psi_{1 p}(z)=p/2+z/4$.
We set the coefficients of these polynomials as
$\phi_{n}(z)=\sum_{m=0}^n a_{n m}z^m$ and
$\psi_{n p}(z)= \sum_{m=0}^n b_{n m}^{(p)}z^m$.
For example,
$$
a_{00}=1, \ a_{1 0}=0, \ a_{1 1}=-\frac{1}{3}, \ a_{2 0}=0, \ 
a_{2 1}=\frac{1}{4}, \ a_{2 2}=\frac{1}{18},
$$
$$
a_{3 0}=0, \ a_{3 1}=-\frac{1}{5}, \ 
a_{3 2}=-\frac{1}{12}, \ a_{3 3}=-\frac{1}{162}, \
b_{0 0}^{(p)}=1, \ b_{1 0}^{(p)}=\frac{p}{2}, \  b_{11}^{(p)}=\frac{1}{4}.
$$
Then the asymptotic behaviors of Hermite polynomials given in \cite{PR29} are summarized as follows.
\vskip 3mm

\noindent {\bf (A.1)} When $x^2 < 2(N+1)$, we put 
$x=\sqrt{2(N+1)}\cos \theta, \ (0<\theta \le \pi/2$). Then for any $L\in \N$
\begin{eqnarray}
&&\frac{H_N(x)}{N!}= \frac{2^{N/2}\exp\big\{(N+1)(1/2 + \cos^2 \theta)\big\}}{(N+1)^{(N+1)/2}(\pi\sin \theta)^{1/2}}
\nonumber\\
&&\times \left[
\sum_{n=0}^{L-1}\sum_{m=0}^n 
C^1_{n m}(N,\theta)\sin \left\{ \frac{N+1}{2}(2\theta-\sin 2\theta) 
+D^1_{n m}(\theta) \right\}
+ {\cal O}\left((N\sin^3\theta)^{-L/2} \right)
\right],
\nonumber
\end{eqnarray}
where
$$
C^1_{n m}(N,\theta)=\frac{1+(-1)^n}{2}
\frac{\Gamma (n + \frac{n+1}{2})}{(N+1)^{n/2}(\sin \theta)^{m+n/2}}a_{n m},
$$
and
$$
D^1_{n m}(\theta)=\frac{\pi}{4} -\frac{\theta}{2} -(2m+n)
\left(\frac{\pi}{4}+\frac{\theta}{2} \right).
$$

\noindent {\bf (A.2)} When $x^2 > 2(N+1)$,  
we put $x=\sqrt{2(N+1)}\cosh \theta, \ (0<\theta < \infty$). 
Then for any $L\in \N$
\begin{eqnarray}
\frac{H_N(x)}{N!}&=& \frac{2^{N/2}\exp\big\{(N+1)(1/2 + \cosh \theta 
(\cosh\theta -\sinh\theta)\big\}}
{(N+1)^{(N+1)/2}(2\pi\sinh \theta)^{1/2}(\cosh\theta -\sinh\theta)^{(N+1/2)}}
\nonumber\\
&\times&\left[
\sum_{n=0}^{L-1}\sum_{m=0}^n C^2_{n m}(\theta, N)
+ {\cal O}\left(N^{-L/2} \theta^{-3L/2})\right)
\right],
\nonumber
\end{eqnarray}
where
$$
C^2_{n m}(\theta, N)= \frac{1+(-1)^n}{2}
\frac{\Gamma (n + \frac{n+1}{2})}{(N+1)^{n/2}}
\left( \frac{-2}{1-e^{-2\theta}} \right)^{m+n/2}a_{n m}.
$$

\noindent {\bf (A.3)} When $x^2 \sim 2(N+1)$, 
we put $x=\sqrt{2(N+1)} -2^{-1/2}N^{-1/6}y$, $y=o(N^{2/3})$. 
Then there exists a positive constant $h^*$ such
\begin{eqnarray}
\frac{H_N(x)}{N!}&=&\frac{e^{3x^2/4}}{\pi \left(x/\sqrt{2}\right)^{N+2/3}}
\left[
\sum_{p=0}^\infty \frac{A_p(x)}{p!}y^p
+{\cal O}\left( x^{-1/3}e^{-h^* x^2}\right)\right]
\nonumber
\end{eqnarray}
with the function $A_p$ having the following asymptotic expansions:
\begin{eqnarray}
&&A_p(x)=3^{(p-2)/3}\sum_{n=0}^{L-1}\sum_{m=0}^n 
\frac{(-1)^m 3^{m+n/3}}{\left( x/\sqrt{2} \right)^{2n/3} }
\Gamma\left(\frac{p+n+1}{3}+m \right)
\sin \left(\frac{p+n+1}{3}\pi\right)b_{n m}^{(p)}
\nonumber\\ 
&&\qquad\qquad\qquad\qquad\qquad\qquad\qquad\qquad\qquad\qquad
+ {\cal O}\left(x^{-2L/3}\right).
\nonumber
\end{eqnarray}
Applying the above to the function 
$\varphi_{N}(x)$ defined in (\ref{def:phin}),
we obtain the lemma. 

\begin{lem}
\label{lemma:Plancherel}
\noindent{\rm (i)} 
For $N\in\N$, $\theta \in [0, \pi/2)$ and $L\in \N$
\begin{eqnarray}
&&\varphi_{N}\Big( \sqrt{2(N+1)} \cos \theta \Big)
= \frac{1+{\cal O}(N^{-1})}{\sqrt{\pi \sin \theta}}
\left( \frac{2}{N} \right)^{1/4}
\nonumber\\
&&\times
\left[
\sum_{n=0}^{L-1}\sum_{m=0}^n C^1_{n m}(N,\theta)
\sin \left\{ \frac{N+1}{2}(2\theta-\sin 2\theta) 
+ D^1_{n m}(\theta)
\right\}
+ {\cal O} \left((N\sin^3\theta )^{-L/2} \right)
\right].
\nonumber
\\
\label{eqn:phiasym1a}   
\end{eqnarray}

\noindent{\rm (ii)}
For $N\in\N$ and $\theta \in [0,\infty)$
\begin{eqnarray}
&& \varphi_{N}\Big( \sqrt{2(N+1)} \cosh \theta \Big)
= \frac{1+{\cal O}(N^{-1})}{\sqrt{2 \pi \sinh \theta}}
\left( \frac{1}{2N} \right)^{1/4} 
\exp \left\{
\left( \frac{N+1}{2} \right)
(2 \theta - \sinh 2 \theta ) \right\}
\nonumber\\
&& \qquad \qquad \qquad \qquad 
\times \left[
\sum_{n=0}^{L-1}\sum_{m=0}^n C^2_{n m}(\theta, N)
+ {\cal O}\Big(N^{-L/2} \theta^{-3L/2})\Big)
\right].
\label{eqn:phiasym1b}
\end{eqnarray}

\noindent{\rm (iii)}
For $N\in\N$ and $y$ with $|y|=o(N^{2/3})$
\begin{eqnarray}
&& \varphi_{N} \left( \sqrt{2(N+1)} - \frac{y}{\sqrt{2} N^{1/6}} \right)
\nonumber\\
&&\qquad\qquad
= 2^{1/4} N^{-1/12} 
\left\{ B \left(y,\sqrt{2(N+1)} 
- \frac{y}{\sqrt{2} N^{1/6}} \right) +{\cal O}(N^{-1}) \right\},
\label{eqn:phiasym1c}
\end{eqnarray}
where
\begin{eqnarray}
B(y,x)& \equiv & {\rm Ai}(-y)+ 
\left(\frac{x}{2}\right)^{-2/3}
\Big\{
c_{10} {\rm Ai}'(-y) + c_{11}y^2{\rm Ai}(-y)
\Big\}
\nonumber\\
&&\qquad+\left(\frac{x}{2}\right)^{-4/3}
\Big\{
c_{20} {\rm Ai}'(-y) + c_{21}y^2{\rm Ai}(-y)
+c_{22}y^3{\rm Ai}'(-y)
\Big\}
\nonumber
\end{eqnarray}
with constants $c_{n m}$, $0\le m \le n \le 2$,
which do not depend on $x$ and $y$.
\end{lem}

\noindent {\it Proof.}
Applying Stirling's formula
$$
N! = \left( \frac{N}{e} \right)^N \sqrt{2\pi N}
\left[ 1+ \frac{1}{12N} + {\cal O}(N^{-2}) \right]
$$
to the results {\bf (A.1)} and {\bf (A.2)}, 
we obtain (\ref{eqn:phiasym1a}) and (\ref{eqn:phiasym1b}) 
by simple calculation.
To obtain (\ref{eqn:phiasym1c}) we use {\bf (A.3)} with $L=3$
and the expansion 
\begin{eqnarray}
{\rm Ai}(-y)&=&\frac{1}{\pi}\sum_{p=0}^\infty 3^{(p-2)/3}
\Gamma \left(\frac{p+1}{3}\right)
\sin \left( \frac{2(p+1)}{3}\pi \right)
\frac{(-y)^p}{p!}
\nonumber\\
&=&\frac{1}{\pi}\sum_{p=0}^\infty 3^{(p-2)/3}
\Gamma \left(\frac{p+1}{3}\right)
\sin \left( \frac{p+1}{3}\pi \right)
\frac{y^p}{p!}
\nonumber
\end{eqnarray}
given in (2.36) in \cite{VS04}.
\qed

\vskip 3mm

We show the main estimate in this section.

\begin{lem}
\label{lemma:keyestimate}
\noindent{\rm (i)} 
Let $x= \sqrt{2N} \cos \theta$ with $N\in\N$ and $\theta \in (0, \pi/2)$.
Suppose that $N\sin^3 \theta \ge CN^{\varepsilon}$
for some $C>0$ and $\varepsilon >0$. Then
\begin{equation}
\sum_{k=0}^{N-1} \Big( \varphi_{k}(x) \Big)^{2} = 
\rho_{\rm sc}^N(1,x)+ {\cal O}\left( \frac{1}{\sqrt{N}\sin^2\theta} \right).
\nonumber
\end{equation}

\noindent{\rm (ii)} 
Let $x= \sqrt{2N} \cosh \theta$ with $N\in\N$ and $\theta >0$.
Suppose that $N\sinh^3 \theta \ge N^{\varepsilon}$ 
for some $C>0$ and $\varepsilon >0$. Then
\begin{equation}
\sum_{k=0}^{N-1} \Big( \varphi_{k}(x) \Big)^{2}= 
{\cal O}\left( \frac{1}{\sqrt{N}\sinh^2 \theta}\right).
\nonumber
\end{equation}

\noindent{\rm (iii)} 
Let $x=2N^{2/3}-y$ with $N\in\N$ and 
$|y|\le CN^{\beta}$ for some $C>0$ and $\beta\in (0,2/21)$
\begin{equation}
\rho_{\rm GUE}^N(N^{1/3},x) = \rho_{\rm Ai}(-y)
+ {\cal O}\left( |y|^{-1} \right).
\nonumber
\end{equation}
\end{lem}

\noindent{\it Proof.} For the proof of this lemma
we use the Christoffel-Darboux formula 
\begin{eqnarray}
\sum_{k=0}^{N-1} \Big( \varphi_{k}(x) \Big)^{2}
&=& N \Big( \varphi_{N}(x) \Big)^2
-\sqrt{N(N+1)} \varphi_{N+1}(x) \varphi_{N-1}(x)
\nonumber\\
&=& N \Big\{ \Big( \varphi_{N}(x) \Big)^2
- \varphi_{N+1}(x) \varphi_{N-1}(x) \Big\}
\left(1+ {\cal O}(N^{-1})\right).
\label{eqn:CD1}
\end{eqnarray}

For proving (i) we first show that for $\ell\in \{ -1, 0, 1\}$
\begin{eqnarray}
&&\varphi_{N+\ell}\Big( \sqrt{2N} \cos \theta \Big)
= \frac{1+{\cal O}(N^{-1})}{\sqrt{\pi \sin \theta}}
\left( \frac{2}{N} \right)^{1/4}
\nonumber\\
&&\qquad\qquad\qquad\times
\Bigg[
\sum_{n=0}^{L-1}\sum_{m=0}^n C^1_{n m}(N-1,\theta)
\sin\bigg\{ \frac{N}{2}(2\theta-\sin 2\theta) 
+ D^1_{n m}(\theta) -(1+\ell)\theta
\bigg\}
\nonumber\\
&&\qquad\qquad\qquad\qquad\qquad\qquad\qquad\qquad\qquad
+ {\cal O}\Big(\frac{1}{N\sin\theta }\Big)\Bigg].
\label{eqn:est3}
\end{eqnarray}
Substituting $N-1$ instead of $N$ in (\ref{eqn:phiasym1a})
and taking $L > 2/\varepsilon$, we have 
(\ref{eqn:est3}) with $\ell=-1$.
For calculating $\varphi_{N+\ell}(x)$ with $\ell\in\{0,1\}$,
we take $\eta_{\ell}$ such that
$$
\cos (\theta+\eta_{\ell})= \sqrt{\frac{2N}{2(N+1+\ell)}} \cos \theta,
$$
and then
\begin{eqnarray}
\varphi_{N+\ell}\Big( \sqrt{2N} \cos \theta \Big)
&= \varphi_{N+\ell} \Big( \sqrt{2(N+1+\ell)} \cos (\theta+\eta_{\ell}) \Big).
\label{eqn:est1}
\end{eqnarray}
Since
$$
\sqrt{\frac{2N}{2(N+1+\ell)}}= \sqrt{1-\frac{1}{N+1+\ell}}
= 1- \frac{1+\ell}{2N} + {\cal O} \left( \frac{1}{N^2} \right)
$$
and
$$
\cos(\theta+\eta_{\ell}) =\cos \theta - \eta_{\ell} \sin \theta
+{\cal O}(\eta_{\ell}^2),
$$
we have
\begin{equation}
\eta_{\ell}= \frac{(1+\ell)\cos \theta}{2N\sin \theta}
+{\cal O} \left(\frac{1}{N^2\sin\theta} \right)
={\cal O} \left(\frac{1}{N\sin\theta} \right)
\label{eqn:eta1}.
\end{equation}
From (\ref{eqn:phiasym1a}) and (\ref{eqn:est1}), we obtain
\begin{eqnarray}
&&\varphi_{N+\ell}\Big( \sqrt{2N} \cos \theta \Big)
= \frac{1+{\cal O}(N^{-1})}{\sqrt{\pi \sin \theta}}
\left( \frac{2}{N} \right)^{1/4}
\nonumber\\
&&\qquad\times
\Bigg[
\sum_{n=0}^{L-1}\sum_{m=0}^n C^1_{n m}(N+\ell,\theta+\eta_{\ell})
\sin\bigg\{ \frac{N}{2}\big(2(\theta+\eta_{\ell})-\sin 2(\theta+\eta_{\ell})\big) 
+ D^1_{n m}(\theta+\eta_{\ell})
\bigg\}
\nonumber\\
&&\qquad\qquad\qquad\qquad\qquad\qquad\qquad\qquad\qquad
+ {\cal O}\Big((N\sin^3\theta )^{-L/2}\Big)
\Bigg].
\label{eqn:est2}
\end{eqnarray}
By simple calculations with (\ref{eqn:eta1}), we have
$$
\frac{N}{2}
\Big\{ \sin 2 (\theta+\eta_{\ell})-2 (\theta+\eta_{\ell}) \Big\}
= \frac{N}{2}(\sin 2 \theta - 2 \theta)
- (1+\ell)\theta 
+{\cal O} \left( \frac{1}{N\sin\theta} \right).
$$
Then from  (\ref{eqn:est2}) with $L > 2/\varepsilon$ and the above estimate,
we obtain (\ref{eqn:est3}) for $\ell\in\{0,1\}$.

Substituting (\ref{eqn:est3}) to (\ref{eqn:CD1})C
from the identity 
$$
2\sin A \sin B - \sin (A+\theta)\sin (B-\theta)-\sin (A-\theta)\sin (B+\theta)
=2\sin^2 \theta \cos (A-B)
$$
with 
$A=N(2\theta-\sin 2\theta)/2+D^1_{n m}(\theta)-\theta$
and
$B=N(2\theta-\sin 2\theta)/2+D^1_{n' m'}(\theta)-\theta$,
$0\le m \le n \le L-1$, $0\le m' \le n' \le L-1$,
we have
\begin{eqnarray}
\sum_{k=0}^{N-1} \Big( \varphi_{k}( \sqrt{2N} \cos \theta ) \Big)^2
&=& \frac{\sqrt{2N}}{\pi}
\sin \theta + {\cal O}\left( \frac{1}{\sqrt{N}\sin^2\theta} \right).
\nonumber
\end{eqnarray}
Since
$\sin \theta = \sqrt{1-\cos^2 \theta}
=\sqrt{1-x^2/\{2(N+1)\}}
= \sqrt{2(N+1)-x^2}/\sqrt{2(N+1)}$,
\begin{equation}
\label{asym:main1}
\sum_{k=0}^{N-1}\Big( \varphi_{k}(x) \Big)^2
= \frac{1}{\pi} \sqrt{2N-x^2}
+ {\cal O}\left( \frac{1}{\sqrt{N}\sin^2\theta} \right)
=\rho_{\rm sc}^N(1,x)+ {\cal O}\left( \frac{1}{\sqrt{N}\sin^2\theta} \right).
\end{equation}
Hence, the proof of (i) is complete.

For proving (ii) we show that for $\ell\in\{-1, 0, 1\}$
\begin{eqnarray}
&& \varphi_{N+\ell} \Big(\sqrt{2N} \cosh \theta \Big)
= \frac{1+{\cal O}(N^{-1})}{\sqrt{2 \pi \sinh \theta}}
\left( \frac{1}{2N} \right)^{1/4} 
\nonumber\\
&& \qquad\qquad\qquad\qquad \times 
\exp 
\left[
\left( \frac{N+1+\ell}{2} \right)
( 2\theta-\sinh 2 \theta) + (1+\ell)\theta \right] 
\nonumber\\
&& \qquad\qquad\qquad\qquad 
\times\Bigg[
\sum_{n=0}^{L-1}\sum_{m=0}^n 
C^2_{n m}(\theta, N+\ell)
+{\cal O}\Big(\frac{\cosh^3\theta}{N\sinh\theta}\Big)
\Bigg].
\label{eqn:est_h2}
\end{eqnarray}
Substituting $N-1$ instead of $N$ in (\ref{eqn:phiasym1b}) 
and taking $L > 2/\varepsilon$, 
we have (\ref{eqn:est_h2}) with $\ell=-1$.
For $\ell\in\{0,1\}$, take $\eta_\ell$ such that
$$
\cosh (\theta+\eta_\ell) = \sqrt{\frac{2N}{2(N+1+\ell)}}\cosh \theta,
$$
and we have
$$
\varphi_{N+1}\Big( \sqrt{2N} \cosh \theta \Big)
= \varphi_{N+1} \Big(
\sqrt{2(N+1+\ell)} \, \cosh (\theta+\eta_\ell) \Big).
$$
Since
$\cosh(\theta+\eta_\ell)=\cosh \theta + \eta_\ell \sinh \theta+{\cal O}(\eta_\ell^2)$,
we have
$$
\eta_\ell=-\frac{(1+\ell)\cosh \theta}{2N\sinh \theta}
+{\cal O} \left( \frac{1}{N^2 \sinh \theta} \right)
={\cal O}\left(\frac{\cosh \theta}{N \sinh \theta}\right).
$$
Substituting the equality to (\ref{eqn:phiasym1b})Cwe have
\begin{eqnarray}
&& \varphi_{N+\ell} \Big(
\sqrt{2N} \cosh \theta \Big)
= \frac{1+{\cal O}(N^{-1})}{\sqrt{2 \pi \sinh \theta}}
\left( \frac{1}{2N} \right)^{1/4} 
\nonumber\\
&& \qquad\qquad\qquad\qquad 
\times \exp \left[
\left( \frac{N+1+\ell}{2} \right)
\Big\{ 2(\theta+\eta_\ell)-\sinh 2 (\theta+\eta_\ell) \Big\}\right] 
\nonumber\\
&& \qquad\qquad\qquad\qquad 
\times\Bigg[
\sum_{n=0}^{L-1}\sum_{m=0}^n C^2_{n m}(\theta+\eta_\ell, N+\ell)
+ {\cal O}\Big(N^{-L/2} \theta^{-3L/2})\Big)
\Bigg].
\nonumber
\end{eqnarray}
Hence if we take $L\in\N$ such that $L>2/\varepsilon$,
from the relation
\begin{eqnarray}
&& 2(\theta+\eta_\ell)-\sinh 2 (\theta+\eta_\ell) 
= 2 \theta - \sinh 2 \theta + \frac{1+\ell}{N} \sinh 2 \theta
+{\cal O}\Big(\frac{\cosh^3\theta}{N^2\sinh\theta}\Big),
\nonumber
\end{eqnarray}
we can conclude (\ref{eqn:est_h2}) with $\ell\in\{0, 1\}$.

Combining  (\ref{eqn:est_h2}) with (\ref{eqn:CD1}),
we obtain 
\begin{eqnarray}
&&\sum_{k=0}^{N-1} \Big(
\varphi_{k} \Big(\sqrt{2N} \cosh \theta \Big) \Big)^2
\nonumber\\
&&= \frac{N}{2 \pi \sinh \theta}
\frac{1}{\sqrt{2N}}
\exp \Big[ (N+1)(2\theta-\sinh 2\theta)+2\theta \Big]
 \times{\cal O} \left( \frac{\cosh^3 \theta}{N\sinh \theta} \right).
\nonumber
\end{eqnarray}
Then (ii) is concluded by simple calculation.


Finally, we show (iii).
Put $w=\sqrt{2N} - 2^{-1/2}N^{-1/6}y$, 
From the Christoffel-Darboux formula (\ref{eqn:CD1})
and (\ref{eqn:phiasym1c}) 
\begin{eqnarray}
&&\varphi_{N-1} \left(\sqrt{2N} - \frac{y}{\sqrt{2} N^{1/6}} \right)
\nonumber\\
&&\qquad= 2^{1/4} N^{-1/12} 
\Big\{ 
B\left( y(1-N^{-1})^{1/6},\sqrt{2N} - \frac{y}{\sqrt{2} N^{1/6}}\right) 
+{\cal O}(N^{-1}) 
\Big\}
\nonumber\\
&&\qquad= 2^{1/4} N^{-1/12} 
\Big\{ 
B\left( y,w \right)
 +{\cal O}(N^{-1}) 
\Big\}.
\nonumber
\end{eqnarray}
Noting the following two simple relations,
$$
\sqrt{2N}-cN^{-1/6}
=\sqrt{2(N+1+\ell)}-\left( c + \frac{1+\ell}{\sqrt{2}N^{1/3}}\right)N^{-1/6}
+{\cal O}(N^{-3/2}),
\quad \ell \in \{0,1\},
$$
and
$y/\sqrt{2} + (1+\ell)/\{\sqrt{2}N^{1/3}\}
= \{y+(1+\ell)N^{-1/3}\}/\sqrt{2}$,
we apply (\ref{eqn:phiasym1c}) to obtain
\begin{eqnarray}
&&\varphi_{N+\ell} \left(\sqrt{2N} - \frac{y}{\sqrt{2} N^{1/6}} \right)
= 2^{1/4} N^{-1/12} 
\left\{ B\left( y+(1+\ell)N^{-1/3}, w- \frac{\ell+1}{\sqrt{2N}}\right) 
+{\cal O}(N^{-1}) \right\}.
\nonumber
\end{eqnarray}
Hence
\begin{eqnarray}
&&\varphi_N
\left( w \right)^2
-\varphi_{N+1} \left( w \right) \varphi_{N-1}\left( w \right)
\nonumber
\\
&&=\sqrt{2}N^{-1/6}
\left[ 
B\bigg(y+N^{-1/3},w-\frac{1}{\sqrt{2N}}\bigg)^2
- B\bigg(y+2N^{-1/3},w-2\frac{1}{\sqrt{2N}}\bigg)B(y,w)
\right]
\nonumber\\
&&\qquad\qquad\qquad\qquad\qquad\qquad\qquad\qquad\qquad\qquad\qquad\qquad\qquad\qquad
+{\cal O}(N^{-7/6}).
\nonumber
\end{eqnarray}
Noting that
\begin{eqnarray}
&&B\bigg(y+N^{-1/3},w-\frac{1}{\sqrt{2N}}\bigg)^2
- B\bigg(y+2N^{-1/3},w-2\frac{1}{\sqrt{2N}}\bigg)B(y,w)
\nonumber\\
&&=\left\{\left(
\frac{\partial }{\partial y}B(y,w)\right)^2
-B(y,w)\frac{\partial^2 }{\partial y^2}B(y,w)
\right\}
N^{-2/3}
\left( 1+ {\cal O}\left(  \frac{y^{1/2}}{wN^{1/6}} \right)\right)
\nonumber
\end{eqnarray}
and
the asymptotic properties of $B(y,x)-{\rm Ai}(-y)$
and its derivatives given by
\begin{eqnarray}
&&B(y,x)-{\rm Ai}(-y)= {\cal O}(x^{-2/3}y^{7/4} + x^{-4/3}y^{13/4})
\nonumber\\
&&\frac{\partial^k}{\partial y^k}\frac{\partial^\ell}{\partial x^\ell}
(B(y,x)-{\rm Ai}(-y))
= {\cal O}(x^{-2/3-\ell}y^{7/4+k/2} + x^{-4/3-\ell}y^{13/4+k/2}),
\nonumber
\end{eqnarray}
$x,y\to\infty$, 
we obtain
\begin{eqnarray}
&&\frac{1}{\sqrt{2} N^{1/6}}\sum_{k=0}^{N-1} 
\varphi_k \left( \sqrt{2N} - \frac{y}{\sqrt{2} N^{1/6}} \right)^2
\nonumber\\
&&=
\left\{
\left(\frac{\partial }{\partial y}B(y,w)\right)^2
-B(y,w)\frac{\partial^2 }{\partial y^2}B(y,w)
\right\}
\left( 1+ {\cal O}\left( \frac{y^{1/2}}{wN^{1/6}} \right)\right)
\nonumber\\
&&=
\left\{
\left(\frac{\partial }{\partial y}{\rm Ai}(-y)\right)^2
-{\rm Ai}(-y)\frac{\partial^2 }{\partial y^2}{\rm Ai}(-y)
+{\cal O}(y^{5/2}w^{-2/3})
\right\}
\left( 1+ {\cal O}\left(\frac{y^{1/2}}{wN^{1/6}}  \right)\right)
\nonumber
\\
&&=
\left(\frac{\partial }{\partial y}{\rm Ai}(-y)\right)^2
-{\rm Ai}(-y)\frac{\partial^2 }{\partial y^2}{\rm Ai}(-y)
+{\cal O}(y^{5/2}N^{-1/3}).
\nonumber
\end{eqnarray}
Hence we conclude that for $x=\sqrt{2}N^{1/6}w=2N^{2/3}-y$, 
$|y|\le N^{\beta}$ with some $\beta\in (0,2/21)$,
\begin{equation}
\rho_{\rm GUE}^N(N^{1/3},x) = 
\left(\frac{\partial }{\partial y}{\rm Ai(-y)}\right)^2
-{\rm Ai}(-y)\frac{\partial^2 }{\partial y^2}{\rm Ai}(-y)
+ {\cal O}\left( |y|^{-1} \right).
\nonumber
\end{equation}
This completes the proof. \qed

\vskip 3mm

Putting $x=2N^{2/3} ( \cos \theta -  1) \in (-\infty,0)$.
we have $-x\sim N^{2/3}\theta^{2}$. 
In case $N\sin^3 \theta \ge N^{3\varepsilon/2}$ for some $\varepsilon >0$, we see from Lemma \ref{lemma:keyestimate} (i) that 
\begin{equation}
\label{estimete:1}
\rho_{\cal A}^N (x)
=\widehat{\rho}_{\rm sc}^N(x) 
+ {\cal O}\left(\frac{1}{N^{2/3}\theta^{2}}\right)
=\widehat{\rho}_{\rm sc}^N(x) 
+ {\cal O}\left(|x|^{-1}\right).
\nonumber
\end{equation}
Similarly, Putting $x=2N^{2/3} ( \cosh \theta -  1) \in (0,\infty)$,
we have $x\sim N^{2/3}(\cosh \theta -1)$. 
In case $N\sinh^3 \theta \ge N^{3\varepsilon/2}$ for some $\varepsilon >0$, we see from Lemma \ref{lemma:keyestimate} (ii) that 
\begin{equation}
\label{estimete:2}
\rho_{\cal A}^N (x)
=\widehat{\rho}_{\rm sc}^N(x) 
+ {\cal O}\left(\frac{1}{N^{2/3}\sinh^2 \theta}\right)
=\widehat{\rho}_{\rm sc}^N(x) 
+ {\cal O}\left(|x|^{-1}\right).
\nonumber
\end{equation}
In the case that $|x|\le N^{\varepsilon}$ for some 
$\varepsilon \in (0, 2/21)$,
we see from Lemma \ref{lemma:keyestimate} (ii) that 
\begin{eqnarray}
\label{estimete:3}
\rho_{\cal A}^N (x)
&=&\rho_{\rm Ai}(x) + {\cal O}\left(|x|^{-1}\right)
\nonumber\\
&=& \widehat{\rho}_{\rm sc}^N(x) 
+ |\widehat{\rho}_{\rm sc}^N(x)-\widehat{\rho}(x)| 
+ |\widehat{\rho}(x)-\rho_{\rm Ai}(x)|
+ {\cal O}\left(|x|^{-1}\right).
\nonumber
\end{eqnarray}
From (\ref{def:circleN}) and (\ref{asym:Airy}) 
we obtain (\ref{asym:AN}).

\vskip 1cm

\begin{small}
\noindent{\bf Acknowledgements} \quad
M.K. is supported in part by
the Grant-in-Aid for Scientific Research (C)
(No.21540397) of Japan Society for
the Promotion of Science.
H.T. is supported in part by
the Grant-in-Aid for Scientific Research 
(KIBAN-C, No.23540122) of Japan Society for
the Promotion of Science.


\end{small}
\end{document}